\documentclass[11pt]{amsart}

\usepackage[all]{xy}
\usepackage{color, pstricks}
\usepackage{mathptmx} 
\usepackage{mathrsfs}
\usepackage{amsmath,amsthm, amssymb, amsgen,amsxtra,amsfonts, amsbsy} 
\usepackage[all]{xy}
\usepackage{verbatim}

\usepackage{times}
\newcommand{\mathds}[1]{{\mathbb #1}}

\begin{document}
%
%
%
\theoremstyle{definition}
\newtheorem{Definition}{Definition}[section]
\newtheorem*{Definitionx}{Definition}
\newtheorem{Convention}{Definition}[section]
\newtheorem{Construction}{Construction}[section]
\newtheorem{Example}[Definition]{Example}
\newtheorem{Examples}[Definition]{Examples}
\newtheorem{Remark}[Definition]{Remark}
\newtheorem*{Remarkx}{Remark}
\newtheorem{Remarks}[Definition]{Remarks}
\newtheorem{Caution}[Definition]{Caution}
\newtheorem{Conjecture}[Definition]{Conjecture}
\newtheorem*{Conjecturex}{Conjecture}
\newtheorem{Question}[Definition]{Question}
\newtheorem{Questions}[Definition]{Questions}
\newtheorem*{Acknowledgements}{Acknowledgements}
\newtheorem*{Organization}{Organization}
\newtheorem*{Disclaimer}{Disclaimer}
\theoremstyle{plain}
\newtheorem{Theorem}[Definition]{Theorem}
\newtheorem*{Theoremx}{Theorem}
\newtheorem{Proposition}[Definition]{Proposition}
\newtheorem*{Propositionx}{Proposition}
\newtheorem{Lemma}[Definition]{Lemma}
\newtheorem{Corollary}[Definition]{Corollary}
\newtheorem*{Corollaryx}{Corollary}
\newtheorem{Fact}[Definition]{Fact}
\newtheorem{Facts}[Definition]{Facts}
\newtheoremstyle{voiditstyle}{3pt}{3pt}{\itshape}{\parindent}%
{\bfseries}{.}{ }{\thmnote{#3}}%
\theoremstyle{voiditstyle}
\newtheorem*{VoidItalic}{}
\newtheoremstyle{voidromstyle}{3pt}{3pt}{\rm}{\parindent}%
{\bfseries}{.}{ }{\thmnote{#3}}%
\theoremstyle{voidromstyle}
\newtheorem*{VoidRoman}{}

%
\newcommand{\prf}{\par\noindent{\sc Proof.}\quad}
\newcommand{\blowup}{\rule[-3mm]{0mm}{0mm}}
\newcommand{\cal}{\mathcal}
\newcommand{\Aff}{{\mathds{A}}}
\newcommand{\BB}{{\mathds{B}}}
\newcommand{\CC}{{\mathds{C}}}
\newcommand{\EE}{{\mathds{E}}}
\newcommand{\FF}{{\mathds{F}}}
\newcommand{\GG}{{\mathds{G}}}
\newcommand{\HH}{{\mathds{H}}}
\newcommand{\NN}{{\mathds{N}}}
\newcommand{\ZZ}{{\mathds{Z}}}
\newcommand{\PP}{{\mathds{P}}}
\newcommand{\QQ}{{\mathds{Q}}}
\newcommand{\RR}{{\mathds{R}}}
\newcommand{\Liea}{{\mathfrak a}}
\newcommand{\Lieb}{{\mathfrak b}}
\newcommand{\Lieg}{{\mathfrak g}}
\newcommand{\Liem}{{\mathfrak m}}
\newcommand{\ideala}{{\mathfrak a}}
\newcommand{\idealb}{{\mathfrak b}}
\newcommand{\idealg}{{\mathfrak g}}
\newcommand{\idealm}{{\mathfrak m}}
\newcommand{\idealp}{{\mathfrak p}}
\newcommand{\idealq}{{\mathfrak q}}
\newcommand{\idealI}{{\cal I}}
\newcommand{\lin}{\sim}
\newcommand{\num}{\equiv}
\newcommand{\dual}{\ast}
\newcommand{\iso}{\cong}
\newcommand{\homeo}{\approx}
\newcommand{\mm}{{\mathfrak m}}
\newcommand{\pp}{{\mathfrak p}}
\newcommand{\qq}{{\mathfrak q}}
\newcommand{\rr}{{\mathfrak r}}
\newcommand{\pP}{{\mathfrak P}}
\newcommand{\qQ}{{\mathfrak Q}}
\newcommand{\rR}{{\mathfrak R}}
%
%
\newcommand{\OO}{{\cal O}}
\newcommand{\numero}{{n$^{\rm o}\:$}}
\newcommand{\mf}[1]{\mathfrak{#1}}
\newcommand{\mc}[1]{\mathcal{#1}}
\newcommand{\into}{{\hookrightarrow}}
\newcommand{\onto}{{\twoheadrightarrow}}
\newcommand{\Spec}{{\rm Spec}\:}
\newcommand{\BigSpec}{{\rm\bf Spec}\:}
\newcommand{\Spf}{{\rm Spf}\:}
\newcommand{\Proj}{{\rm Proj}\:}
\newcommand{\Pic}{{\rm Pic }}
\newcommand{\Br}{{\rm Br}}
\newcommand{\NS}{{\rm NS}}
\newcommand{\Sym}{{\mathfrak S}}
\newcommand{\Aut}{{\rm Aut}}
\newcommand{\Autp}{{\rm Aut}^p}
\newcommand{\Hom}{{\rm Hom}}
\newcommand{\Ext}{{\rm Ext}}
\newcommand{\ord}{{\rm ord}}
\newcommand{\coker}{{\rm coker}\,}
\newcommand{\divisor}{{\rm div}}
\newcommand{\Def}{{\rm Def}}
\newcommand{\piet}{{\pi_1^{\rm \acute{e}t}}}
\newcommand{\Het}[1]{{H_{\rm \acute{e}t}^{{#1}}}}
\newcommand{\Hfl}[1]{{H_{\rm fl}^{{#1}}}}
\newcommand{\Hcris}[1]{{H_{\rm cris}^{{#1}}}}
\newcommand{\HdR}[1]{{H_{\rm dR}^{{#1}}}}
\newcommand{\hdR}[1]{{h_{\rm dR}^{{#1}}}}
\newcommand{\defin}[1]{{\bf #1}}
\newcommand{\oX}{\cal{X}}
\newcommand{\oA}{\cal{A}}
\newcommand{\oY}{\cal{Y}}

\title{Supersingular K3 Surfaces Are Unirational}
\author{Christian Liedtke}
\address{TU M\"unchen, Zentrum Mathematik - M11, Boltzmannstr. 3, D-85748 Garching bei M\"unchen, Germany}
\curraddr{}
\email{liedtke@ma.tum.de}

\date{April 8, 2014}

\begin{abstract}
   We show that supersingular K3 surfaces 
   in characteristic $p\geq5$ 
   are related by purely inseparable isogenies.
   This implies that they are unirational,
   which proves conjectures of Artin, Rudakov, Shafarevich, and Shioda.
   As a byproduct, we exhibit the moduli space of rigidified K3 crystals
   as an iterated $\PP^1$-bundle over $\FF_{p^2}$.
   To complete the picture, we also establish Shioda--Inose type isogeny theorems 
   for K3 surfaces with
   Picard rank $\rho\geq19$ in positive characteristic.
\end{abstract}

\maketitle

\section{Introduction}
\footnote{Due to a mistake in Proposition 3.5, our results merely show that supersingular K3 surfaces are related by sequences of very special correspondences, which is not sufficient to deduce their unirationality. The unirationality conjecture thus remains a conjecture. See the appendix for erratum.}
The Picard rank $\rho$ of a complex K3 surface satisfies $\rho\leq20$.
In \cite{shioda inose}, \cite{inose}, Shioda and Inose classified complex
K3 surfaces with Picard rank $20$, so-called {\em singular K3 surfaces},
showed that they can be defined over numbers, and thus, form a countable
set and have no moduli.
They also showed that such a surface rationally dominates and 
is rationally dominated  by a Kummer surface.
This is related to a conjecture of Shafarevich \cite{shafarevich},
according to which every 
Hodge-isogeny between the transcendental lattices of two complex K3 surfaces is 
induced by a rational map or a rational correspondence --  we refer to
Section \ref{subsec: isogeny} for details.

The first result of this article
is an extension of the Shioda--Inose theorem to positive characteristic:

\begin{Theoremx}
  Let $X$ be a K3 surface in odd characteristic with Picard rank $19$ or $20$.
  Then, there exists an ordinary Abelian surface $A$ and dominant, rational maps
  $$
      {\rm Km}(A)\,\dashrightarrow\,X\,\dashrightarrow\,{\rm Km}(A),
  $$
  both of which are generically finite of degree $2$.
\end{Theoremx}

We refer to Theorem \ref{shioda inose in p} for more precise statements, fields
of definition, as well as lifting results.
For example, singular K3 surfaces in odd characteristic can be defined over finite fields,
and thus, also these surfaces form a countable set and have no moduli.

Artin \cite{artin} observed that there do not exist K3 surfaces with 
Picard rank $21$ in any characteristic.
On the other hand, Tate \cite{tate} and Shioda \cite{shioda some results}
gave examples of K3 surfaces with Picard rank $22$ in positive
characteristic, so-called {\em Shioda-supersingular K3 surfaces}.
Artin \cite{artin} showed that Shioda-supersingular K3 surfaces are
{\em Artin-supersingular}, that is, their formal Brauer groups are of infinite
height.
It follows from recent progress in the Tate-conjecture for K3 surfaces 
due to Charles \cite{charles}, Madapusi Pera \cite{madapusi}, and 
Maulik \cite{maulik} that a K3 surface in odd characteristic is Artin-supersingular
if and only if it is Shioda-supersingular.

Artin \cite{artin} also showed that supersingular K3 surfaces form $9$-dimensional families, 
which is in contrast to the above mentioned rigidity of singular K3 surfaces.
Moreover, Shioda \cite{shioda some results} showed that Tate's and his  examples 
are unirational, another property of K3 surfaces that can happen in positive characteristic only.
Since unirational K3 surfaces are supersingular as shown by Shioda \cite{shioda an example}, 
this led several people to conjecture the converse:

\begin{Conjecturex}[Artin, Rudakov, Shafarevich, Shioda]
  A K3 surface is supersingular if and only if it is unirational.
\end{Conjecturex}

Shioda \cite{shioda some results} established this conjecture 
for supersingular Kummer surfaces in odd characteristic,  
Rudakov and Sharafevich \cite{rudakov shafarevich} showed it 
in characteristic $2$ and for K3 surfaces with Artin invariant  $\sigma_0\leq6$ in characteristic $3$, 
and Pho and Shimada \cite{shimada pho} for K3 surfaces with
Artin invariant $\sigma_0\leq3$ in characteristic $5$.
We refer to \cite{katsura} and \cite{ito liedtke} for some refinements.
In particular, there do exist unirational K3 surfaces
in every positive characteristic.
\medskip

The key result of this article is a structure theorem for supersingular K3
surfaces, which was posed as a question by 
Rudakov and Shafarevich in \cite{rudakov shafarevich}, and
which is similar to the Shioda--Inose theorem for singular K3 surfaces.

\begin{Theoremx}
  Let $X$ and $X'$ be supersingular K3 surfaces in characteristic $p\geq5$
  with Artin invariants $\sigma_0$ and $\sigma_0'$, respectively.
  Then, there exist dominant and rational maps
      $$
        X \,\dashrightarrow\,X'\,\dashrightarrow\,X\,,
      $$
      both of 
      which are purely inseparable and generically finite of degree $p^{2\sigma_0+2\sigma_0'-4}$.
\end{Theoremx}

In \cite{shioda inose}, \cite{inose}, Shioda and Inose introduced a notion of
{\em isogeny} for singular K3 surfaces over the complex numbers, which was extended
to other types of complex K3 surfaces by Morrison \cite{morrison},
Mukai \cite{mukai}, and Nikulin \cite{nikulin rational}.
We refer to Section \ref{subsec: isogeny} for an extension of this notion to positive characteristic,
and using this terminology, our structure theorem says that all supersingular K3 
surfaces are mutually {\em purely inseparably isogenous}.

Our theorem also fits into Shafarevich's conjecture  \cite{shafarevich}
mentioned above:
supersingular K3 surfaces are precisely those K3 surfaces without transcendental 
cycles in their second $\ell$-adic cohomology. 
Thus, their ``transcendental lattices'' should be thought of 
as being zero, thus mutually isogenous, and by our theorem, they are 
all related by rational maps.
We refer to Section \ref{subsec: isogeny} for details.

Our theorem also explains why supersingular K3 surfaces form $9$-dimensional families, whereas
singular K3 surface have no moduli: in both cases, these surfaces are isogenous to Kummer surfaces.
For singular K3 surfaces, the isogeny is separable and does not deform.
For supersingular K3 surfaces, the isogeny can be chosen purely inseparable,
and deforms in families.
We refer to Remark \ref{rem:correspondences} for details.

The main tool to proving this theorem is that a Jacobian elliptic fibration 
$X\to\PP^1$ on a supersingular
K3 surface with Artin invariant $\sigma_0$ admits a deformation that is
a one-dimensional family of elliptic supersingular K3 surfaces, such that all 
elliptic fibrations in this family are generically torsors under $X\to\PP^1$.
We call this a {\em moving torsor family}  and
refer to Section \ref{subsec:moving torsor general} for details.
Moreover, the generic fiber of this family has Artin invariant $\sigma_0+1$ and is related to
the special fiber $X$ by a purely inseparable isogeny,
see Theorem \ref{thm:moving torsors general}.

In \cite{ogus}, Ogus introduced moduli spaces ${\cal M}_N$ 
of $N$-rigidified K3-crystals, where $N$ is a supersingular K3 lattice.
If $N$ and $N_+$ denote supersingular K3 lattices in odd characteristic
of Artin invariants $\sigma_0$ and $\sigma_0+1$, respectively,
then these moving torsor families induce a structure of a $\PP^1$-bundle.

\begin{Theoremx}
  There exists a surjective morphism 
  $$
    {\cal M}_{N_+}\,\to\,{\cal M}_N
  $$
  together with a section, which turns ${\cal M}_{N_+}$ into a $\PP^1$-bundle over ${\cal M}_N$.
\end{Theoremx}

Using Ogus' Torelli theorem \cite{ogus torelli}, we use this $\PP^1$-bundle structure to show 
that every supersingular K3 surface of
Artin invariant $\sigma_0+1$ is purely inseparably isogenous to one of Artin invariant $\sigma_0$,
and, by induction on the Artin invariant, we obtain our theorem. 
We refer to  Theorem \ref{thm: bundle structure} and Theorem \ref{thm: moduli interpretation} 
for details.
\medskip

As already mentioned above, Shioda \cite{shioda some results} proved that
supersingular Kummer surfaces in odd characteristic are unirational.
Combined with our structure theorem, 
this establishes the Artin--Rudakov--Shafarevich--Shioda conjecture.

\begin{Theoremx}
  Supersingular K3 surfaces in characteristic $p\geq5$ are unirational.
\end{Theoremx}

Together with results of Artin, Shioda, and the recent proof of the
Tate-conjecture for K3 surfaces in odd characteristic,  
we obtain the following equivalence.

\begin{Theoremx}
 For a K3 surface $X$ in characteristic $p\geq5$, the following conditions
 are equivalent:
 \begin{enumerate}
  \item $X$ is unirational.
  \item The Picard rank of $X$ is $22$.
  \item The formal Brauer group of $X$ is of infinite height.
  \item For all $i$, the $F$-crystal $\Hcris{i}(X/W)$ is of slope $i/2$.
 \end{enumerate}
\end{Theoremx}

We refer to Section \ref{subsec: small char}, Section \ref{subsec: ogus small char},
and Section \ref{subsec: small char 2}  for partial results in small characteristics.
For example, if the
Rudakov--Shafarevich theorem \cite{rudakov shafarevich degeneration} on 
potential good reduction of supersingular K3 surfaces were known to hold 
in characteristic $3$, then the above theorems would hold in characteristic $3$
as well.
\medskip

Recently, Lieblich \cite{lieblich unirational k3} gave a new proof of the unirationality of supersingular 
K3 surfaces, following a simplified form of our strategy .
\medskip

This article is organized as follows:

In Section \ref{sec:shioda-inose}, after reviewing formal Brauer groups, several
notions of supersingularity,
and introducing purely inseparable isogenies, 
we classify K3 surfaces with Picard ranks $19$ and $20$ 
in odd characteristic, which generalizes the classical Shioda--Inose theorem.

In Section \ref{sec:torsors}, we show how a supersingular K3 surface with Artin invariant
$\sigma_0$ together with a Jacobian elliptic fibration gives rise to a one-dimensional family 
of elliptic supersingular K3 surfaces that are generically torsors under this Jacobian fibration,
and whose generic fiber has Artin invariant $\sigma_0+1$.
Moreover, we show how these torsors are related to the trivial torsor 
by purely inseparable isogenies.

In Section \ref{sec: ogus moduli}, we interpret these one-dimensional families of torsors
in  terms of Ogus' moduli spaces of supersingular K3 crystals.
As an interesting byproduct, we find that these moduli spaces are related
to each other by iterated $\PP^1$-bundles, together with a moduli interpretation 
of this structure. In particular, this gives a new description of these moduli spaces.

In Section \ref{sec:shioda}, we use the results of the previous sections to prove
that supersingular K3 surfaces are connected by purely inseparable isogenies.
Since Shioda showed that supersingular Kummer surfaces are unirational,
we conclude that all supersingular K3 surfaces are unirational.
Finally, we also characterize unirational Enriques surfaces.

\begin{Acknowledgements}
 It is a pleasure for me to thank 
 Xi Chen, Igor Dolgachev, Gerard van der Geer, 
 Brendan Hassett, Daniel Huybrechts, Toshiyuki Katsura, Frans Oort, Matthias Sch\"utt, 
 Tetsuji Shioda, Burt Totaro, and the referee for discussions, comments,
 and pointing out inaccuracies.
 I especially thank Olivier Benoist and Max Lieblich for pointing out mistakes 
 in earlier versions of this article and helping me to fix them.
 I gratefully acknowledge funding from DFG via Transregio SFB 45, as part of this
 article was written while staying at Bonn university.
\end{Acknowledgements}

\section{Non-supersingular K3 surfaces with large Picard number}
\label{sec:shioda-inose}

In this section, we first review the formal Brauer group, and discuss several notions
of supersingularity for K3 surfaces.
Then, we classify non-supersingular K3 surfaces with large Picard rank
in positive characteristic, which establishes a structure result similar to the Shioda--Inose
theorem over the complex numbers.

\subsection{Formal Brauer groups, supersingularity, and Picard ranks}
\label{subsec: brauer}
Let $X$ be a K3 surface over a field $k$.
By results of Artin and Mazur \cite{artin mazur}, 
the functor on local Artinian $k$-algebras with residue field $k$ 
defined by
$$
 \begin{array}{ccccc}
  \Phi_{X/k}^2 &:& (\mbox{Art}/k) & \to & (\mbox{Abelian groups})\\
     & & R &\mapsto&\ker\left( \Het{2}(X\times_k\Spec R, \GG_m)\,\to\, \Het{2}(X,\GG_m) \right)
 \end{array}
 $$
is pro-representable by a one-dimensional formal group law $\widehat{{\rm Br}}(X)$, which 
is called the {\em formal Brauer group}.
Over algebraically closed fields of positive characteristic,
one-dimensional formal group laws are classified by their {\em height}, and
Artin \cite[Theorem (0.1)]{artin} showed that
the height $h$ of the formal Brauer group of a K3 surface satisfies
$1\leq h\leq10$ or $h=\infty$.

\begin{Definition}
 \label{def: artin-supersingular}
 Let $X$ be a K3 surface over a field of positive characteristic
 and let $h$ be the height of its formal Brauer group.
 Then, $X$ is called {\em ordinary} if $h=1$, and $X$ is called {\em Artin-supersingular}
 if $h=\infty$.
\end{Definition}

The general picture is as follows: 
a smooth and proper variety $X$ over a perfect field  of positive characteristic
is called {\em ordinary}, if the Hodge- and Newton-polygons of all its crystalline
cohomology groups coincide.
It is called {\em supersingular}, if the Newton-polygons of all its crystalline cohomology
groups are straight lines, that is, if the $F$-crystal $\Hcris{i}(X/W)$ is of
slope $i/2$ for all $i$.
Now, if $X$ is a K3 surface, this general definition translates into a condition on $\Hcris{2}(X/W)$
only.
More precisely, being ordinary translates into having slopes $(0,1,2)$, 
and being supersingular into being of slope $1$.
By \cite[Section II.7.2]{illusie}, the slopes of $\Hcris{2}(X/W)$ 
in terms of the height $h$
of the formal Brauer group are $(1-\frac{1}{h},1,1+\frac{1}{h})$.
Thus, for K3 surfaces, Definition \ref{def: artin-supersingular} coincides
with the general picture.

For surfaces, Shioda \cite{shioda an example} introduced another notion
of supersingularity.
To explain it, we note that the 
first Chern class map $c_1:\NS(X)\to H^1(\Omega_X^1)$
is injective over the complex numbers, which implies that 
the Picard rank $\rho$ of a smooth complex projective variety is bounded
above by $h^1(\Omega_X^1)$.
For complex K3 surfaces, this gives the estimate
$\rho\leq20$.
In positive characteristic, Igusa \cite{igusa} established the 
inequality $\rho\leq b_2$, which, for K3 surfaces, only gives the 
estimate $\rho\leq22$.
However, this bound is sharp, since Tate \cite{tate} and Shioda \cite{shioda some results}
showed that there do exist K3 surfaces  
with Picard rank $22$ in positive characteristic.

\begin{Definition}
  \label{def: shioda supersingular}
  Let $X$ be a K3 surface over an algebraically closed field.
  Then, $X$ is called {\em singular} if $\rho=20$, and it is called
  {\em Shioda-supersingular} if $\rho=22$.
\end{Definition}

The relation between these two notions of supersingularity is as follows:
In \cite[Theorem (0.1)]{artin}, Artin showed that a K3 surface whose 
formal Brauer group is of finite height $h$ satisfies $\rho\leq b_2-2h$.
This implies that Shioda-supersingular K3 surfaces are Artin-supersingular.
In \cite[Theorem (4.3)]{artin}, Artin proved that 
Artin-supersingular K3 surfaces that are elliptic are Shioda-supersingular.
In general, the equivalence of Artin- and Shioda-supersingularity
follows from the Tate-conjecture for supersingular K3 surfaces.
Since this has been recently established in odd characteristic by
Charles \cite{charles}, Madapusi Pera \cite{madapusi}, and Maulik \cite{maulik},
we can summarize these results as follows.

\begin{Theorem}[Artin, Charles, Madapusi Pera, Maulik, et al.]
 \label{thm: equivalence}
 For a K3 surface $X$ in odd characteristic, the following are equivalent:
 \begin{enumerate}
 \item $X$ is Shioda-supersingular, that is, $\rho=22$.
 \item $X$ is Artin-supersingular, that is, $h(\widehat{{\rm Br}}(X))=\infty$.
  \item For all $i$, the $F$-crystal $\Hcris{i}(X/W)$ is of slope $i/2$.\qed
 \end{enumerate}
\end{Theorem}

By \cite[Section 4]{artin}, the 
discriminant of the N\'eron--Severi lattice of a Shioda-supersingular K3 surface
is equal to $-p^{2\sigma_0}$ for some integer $1\leq \sigma_0\leq 10$.

\begin{Definition}
  The integer $\sigma_0$ is called the {\em Artin-invariant} of $X$.
\end{Definition}

The Artin invariant $\sigma_0$ gives rise to a stratification of the moduli
space of Shioda-supersingular K3 surfaces \cite[Section 7]{artin}, and it determines the N\'eron--Severi
lattice of a Shioda-supersingular K3 surface up to isometry \cite[Section 1]{rudakov shafarevich}.
We refer the interested reader to the overview articles by Shioda \cite{shioda overview} and 
Rudakov--Shafarevich \cite{rudakov shafarevich overview} for basic properties of 
Shioda-supersingular K3 surfaces, details and further references.

\subsection{Isogenies between K3 surfaces}
\label{subsec: isogeny}
For Abelian varieties, the notion of isogeny  is classical.
For K3 surfaces, there are several and conflicting extensions
of this notion, and we refer to \cite[Section 1]{morrison isogeny}
for an overview.
Following Inose \cite{inose}, we use the most naive one, which
is sufficient for the purposes of this article.

\begin{Definition}
  \label{def: isogeny}
  Let $X$ and $Y$ be varieties of the same dimension
  over a perfect field of positive characteristic $p$.
  An {\em isogeny of degree $n$} from $X$ to $Y$ is a dominant, rational, and generically finite 
  map $X\dashrightarrow Y$ of degree $n$. 
  A {\em purely inseparable isogeny of height $h$} is an isogeny that is purely inseparable 
  of degree $p^h$.
\end{Definition}

For Abelian varieties $A$, $B$ and an isogeny $A\to B$,
there exists an integer $n$ such that multiplication by 
$n:A\to A$ factors through this isogeny.
Such a factorization gives rise to an isogeny $B\to A$, and 
in particular, being isogenous is an equivalence relation.
Over the complex numbers, K3 surfaces with Picard rank $20$
are related to Kummer surfaces by isogenies, and the 
existence of an isogeny in the other direction is a true, but 
non-trivial fact, see \cite{shioda inose}, \cite{inose}, and
\cite{ma}.

Coming back to Definition \ref{def: isogeny}, if
$X\dashrightarrow Y$ is a purely inseparable isogeny of height $h$,
the $h$-fold $k$-linear Frobenius $F^h:X\to X^{(p)}$ factors through this isogeny, 
inducing an isogeny $Y\dashrightarrow X^{(p)}$,
which is purely inseparable of height $(d-1)h$, where $d$
is the dimension of $X$ and $Y$. 
Since the ground field is perfect, we can identify
$X$ with $X^{(p)}$, and thus,
purely inseparable isogenies define an equivalence relation.

Since it motivates some of our results later on and sheds another light
on them, let us shortly discuss a conjecture of Shafarevich concerning complex
K3 surfaces:
let $X$ and $Y$ be complex K3 surfaces with transcendental lattices $T(X)$
and $T(Y)$.
If $\rho(X)=\rho(Y)=20$, then $T(X)$ and $T(Y)$ are of rank $2$, and 
the Shioda--Inose theorem \cite{shioda inose}
says that every isogeny $T(X)\to T(Y)$ preserving Hodge structures
induces and is induced by an isogeny between 
the corresponding surfaces.
Morrison \cite{morrison},
Mukai \cite{mukai}, and Nikulin
\cite{nikulin russian}, \cite{nikulin rational} generalized these results
to K3 surfaces, whose transcendental lattices are of higher rank.
Moreover, Shafarevich \cite{shafarevich} conjectured that every
Hodge isogeny between 
transcendental lattices of complex K3 surfaces is induced
by an isogeny, or, by a rational correspondence.
Here, the right definition of isogeny for K3 surfaces is one difficulty, and
we refer to \cite[Section 1]{morrison isogeny} for a discussion and
the relation of Shafarevich's conjecture to the Hodge conjecture.
We note that results of Chen \cite{chen} imply that 
Shafarevich's conjecture cannot be true if one only allows  
isogenies in the sense of our naive Definition \ref{def: isogeny}.

In positive characteristic, a K3 surface $X$ is Shioda-supersingular if and only
if every class in $\Het{2}(X,\QQ_\ell)$ is algebraic if and only if
the cokernel of $c_1:\NS(X)\to\Hcris{2}(X/W)$ is a $W$-module that is torsion.
Therefore, the ``transcendental lattices'' of Shioda-supersingular K3 surfaces
should be thought of as being zero, in which case they would all be isogenous
for trivial reasons.
Now, if one believes in a characteristic-$p$ version
of Shafarevich's conjecture (even whose precise formulation is unclear
to the author at the moment), one might expect that all Shioda-supersingular
K3 surfaces are related by isogenies.
This was posed as Question 8 by Rudakov and Shafarevich
at the end of \cite{rudakov shafarevich}, and we shall prove it in
Theorem \ref{thm:correspondences} below.

\subsection{The Shioda--Inose theorem in odd characteristic}
In this subsection, 
we classify non-supersingular K3 surfaces
with Picard rank $\rho\geq19$ in odd characteristic and establish
an analog of the Shioda--Inose theorem \cite{shioda inose}, \cite{inose}.
The idea  is to show that such surfaces are ordinary, 
which implies that we can lift them to the Witt ring together with their
Picard groups.
Then, we use the Shioda--Inose theorem in characteristic zero to deduce 
the structure result in odd characteristic.

\begin{Theorem}
  \label{shioda inose in p}
  Let $X$ be a K3 surface with Picard rank $19\leq\rho\leq21$ over an algebraically closed field $k$
  of odd characteristic.
  Then,
  \begin{enumerate}
   \item  $X$ is an ordinary K3 surface, and
   \item  $X$ lifts projectively together with its Picard group to $\Spec W(k)$. 
  \end{enumerate}
  Moreover, 
  \begin{enumerate}
   \setcounter{enumi}{2}
   \item If $\rho=19$, then there exists an ordinary Abelian surface $A$ over $k$, and
    isogenies of degree $2$ 
    $$
      {\rm Km}(A) \,\dashrightarrow\, X \,\dashrightarrow\, {\rm Km}(A)\,.
    $$
   Moreover, neither $X$ nor $A$ can be defined over a finite field.
   \item If $\rho=20$, then there exist two ordinary and isogenous elliptic curves $E$ and $E'$ over $k$,
    and isogenies of degree $2$
    $$
      {\rm Km}(E\times E') \,\dashrightarrow\, X \,\dashrightarrow\, {\rm Km}(E\times E')\,.
    $$
   Moreover, $X$ can be defined over a finite field.
   The lift of $(X,\Pic(X))$ is unique and coincides with the canonical Serre--Tate lift of $X$.
  \item K3 surfaces with Picard rank $\rho=21$ do not exist.
 \end{enumerate}
\end{Theorem}

\begin{Remark}
  Non-existence of K3 surfaces with Picard rank $21$ was already observed by Artin \cite[p. 544]{artin}.
  Independently, Jang \cite[Section 4]{jang} obtained a similar classification result for K3 surfaces with
  Picard rank $\rho=20$.
\end{Remark}

\prf
We proceed in several steps:
\medskip

{\sc Step 1}: $X$ is ordinary and $\rho=21$ is impossible.

Let $h$ be the height of the formal Brauer group.
Since $\rho<22$, we deduce $h<\infty$ from \cite[Theorem 1.7]{artin}.
But then, the inequalities $\rho\leq b_2-2h\leq20$ from \cite[Theorem 0.1]{artin} show
that $\rho=21$ is impossible.
They also show that $X$ is ordinary, that is, $h=1$, if $19\leq\rho\leq20$.
This establishes claims (1) and (5).
\medskip

{\sc Step 2}: There exists a projective lift of the pair $(X,\Pic(X))$ to $W(k)$.

Since $X$ is ordinary, there exists a canonical formal lift ${\cal X}\to\Spf W(k)$, the Serre--Tate lift.
By \cite[Proposition 1.8]{nygaard}, it has the property that $\Pic(X)$ lifts to $\cal X$.
In particular, lifting an ample invertible sheaf, it follows from Grothendieck's existence
theorem 
that ${\cal X}$ is algebraizable.
This establishes claim (2).
\medskip

{\sc Step 3}: If $\rho=20$, then ${\cal X}_{\overline{K}}$ is dominated by a Kummer surface.

Let us now assume that $\rho=20$, let $K$ be the field of fractions of $W(k)$, and 
let $\overline{K}$ be its algebraic closure.
By the classical Shioda--Inose theorem from \cite{shioda inose} and \cite{inose},
there exist isogenous elliptic curves $\widetilde{E}$ and $\widetilde{E}'$ 
with complex multiplication
over $\overline{K}$, and a symplectic involution $\imath$ on
the Kummer surface ${\rm Km}(\widetilde{E}\times \widetilde{E}')$, such that
${\cal X}_{\overline{K}}$ is the desingularization of 
the quotient ${\rm Km}(\widetilde{E}\times \widetilde{E}')/\langle\imath\rangle$.
\medskip

{\sc Step 4}: This Kummer surface has potential good reduction and $\imath$ extends.

Since elliptic curves with complex multiplication have potential good reduction,
there exists a model of ${\rm Km}(\widetilde{E}\times \widetilde{E}')$ over 
a finite extension $R\supseteq W(k)$
with good reduction that is itself a Kummer surface, say, 
${\rm Km}({\cal E}\times{\cal E}')$
(since $p\neq2$, we can form the quotient by the sign involution
over $R$ without trouble).
After possibly enlarging $R$, the involution $\imath$ is defined on the generic fiber 
${\rm Km}({\cal E}\times{\cal E}')_K$.
Now, $\imath$ extends to an involution on ${\rm Km}({\cal E}\times{\cal E}')$, see, 
for example the proof of \cite[Theorem 2.1]{lieblich}.
Since $\imath$ acts trivially on the global $2$-form of the generic fiber, 
its extension will act trivially on the global $2$-form of the special fiber,
and thus, $\imath$ extends to a symplectic involution on 
${\rm Km}({\cal E}\times{\cal E}')\to\Spec R$.
On the geometric generic fiber 
it has precisely $8$ fixed points by \cite{nikulin} or \cite[Lemma 5.2]{morrison}, 
and the same is true for the induced involution on the special fiber 
by \cite[Theorem 3.3]{dolgachev} (here, we use again that $p\neq2$).
\medskip

{\sc Step 5}: $X$ is the quotient of a Kummer surface by an involution.

After possibly enlarging $R$ again, we may form the quotient 
${\rm Km}({\cal E}\times{\cal E}')/\langle \imath\rangle$ and resolve the resulting $8$ families of 
$A_1$-singularities to obtain a smooth family ${\cal Y}\to\Spec R$.
After possibly enlarging $R$ again, the generic fibers of $\cal X$ and $\cal Y$ 
become isomorphic.
Since ${\cal X}$ and ${\cal Y}$ both have good reduction, and
their special fibers are not ruled, the special fibers are isomorphic by the 
Matsusaka--Mumford theorem \cite[Theorem 2]{matsusaka}.
This shows the existence of a rational dominant map 
${\rm Km}(E\times E')\dashrightarrow X$, which is generically finite of degree $2$.
Here, $E$ and $E'$ denote the reductions of $\cal E$ and ${\cal E}'$,
respectively.
The existence of a rational dominant map $X\dashrightarrow {\rm Km}(E\times E')$,
generically finite of degree $2$, follows from the corresponding characteristic zero statement 
as before and we leave the proof to the reader.
\medskip

{\sc Step 6}: Ordinarity and fields of definition.

Since $X$ is ordinary, Frobenius acts bijectively on $H^2(X,\OO_X)$, from which we conclude
that it also acts bijectively on $H^2({\rm Km}(E\times E'), \OO_{{\rm Km}(E\times E')})$ and thus,
on $H^2(E\times E',\OO_{E\times E'})$.
In particular, $E\times E'$ is an ordinary Abelian surface, which implies that
$E$ and $E'$ are ordinary elliptic curves.
And finally, since $\widetilde{E}$ and $\widetilde{E}'$ are elliptic curves with complex multiplication,
they can be defined over $\overline{\QQ}$, which implies that 
${\cal E}$, ${\cal E}'$, ${\rm Km}({\cal E}\times{\cal E}')$  and $\imath$ can be defined over $W(\overline{\FF}_p)$,
which implies that $E$, $E'$, ${\rm Km}(E\times E')$ and $X$ can be defined over $\overline{\FF}_p$.
This establishes claim (4).
\medskip

{\sc Step 7}: Sketch of the case of Picard rank $\rho=19$.

As in step 2, let ${\cal X}\to\Spec W(k)$ be a projective lift of $(X,\Pic(X))$.
Then, as in step 3, there exists an Abelian variety $\widetilde{A}$ over some 
finite extension $L\supseteq K$ and an involution $\imath$ on ${\rm Km}(\widetilde{A})$ such that
${\rm Km}(\widetilde{A})/\imath$ and ${\cal X}_K$ become isomorphic over $\overline{K}$.
Since ${\cal X}$ has good reduction,
the Galois-action of  $G_L:={\rm Gal}(\overline{K}/L)$
on $\Het{2}({\cal X}_{\overline{K}},\QQ_\ell)$, $\ell\neq p$, is unramified.
From this, it is not difficult to see that also the $G_L$-actions on 
$\Het{2}({\rm Km}(\widetilde{A})_{\overline{K}},\QQ_\ell)$ and 
$\Het{2}(\widetilde{A}_{\overline{K}},\QQ_\ell)$ are unramified.
Thus, by the N\'eron--Ogg--Shafarevich criterion,
there exists a smooth model of $\widetilde{A}$ over some finite extension of $W(k)$, 
whose special fiber $A$ is an Abelian surface.
As in step 5, we find rational dominant maps
${\rm Km}(A)\dashrightarrow X$ and $X\dashrightarrow{\rm Km}(A)$, both
of which are generically finite of degree $2$.
As in step 6, we conclude that $A$ is an ordinary Abelian surface.
Finally, if $X$ were definable over $\overline{\FF}_p$, then its geometric Picard rank 
would be even by \cite[p. 544]{artin}, a contradiction. 
This implies that $A$, ${\rm Km}(A)$, and $X$
cannot be defined over $\overline{\FF}_p$
and establishes claim (3).
\qed\medskip

\begin{Remark}
 We would like to point out the following analogy between zero and positive
 characteristic for K3 surfaces with Picard rank $20$:
 over the complex numbers, such surfaces can be defined over $\overline{\QQ}$,
 and thus, form a countable set and have no moduli.
 In characteristic $p\geq3$, such surfaces can be defined
 over $\overline{\FF}_p$, and again, form a countable set and have no moduli.
\end{Remark}

\section{Continuous Families of Torsors}
\label{sec:torsors}

In this section, we consider Jacobian (quasi-)elliptic
fibrations on surfaces in positive characteristic $p$.
If the formal Brauer group of the surface is not $p$-divisible,
then we construct a deformation of the 
Jacobian to a non-Jacobian fibration,
which is generically a family of torsors under the Jacobian fibration.
Using a purely inseparable multisection, we show that 
the special and the generic fiber of this family are related 
by a purely inseparable isogeny.
Our main result is Theorem \ref{thm:moving torsors general},
which is the technical heart of this article.
For a K3 surface, such a family exists if only if it is supersingular
with Artin invariant $\sigma_0\leq9$, and then,
this family can be spread out to a 
smooth family of supersingular K3 surfaces over a proper curve 
such that the generic fiber has Artin invariant $\sigma_0+1$.

In order to avoid confusion, let us fix the following terminology.

\begin{Definition}
 A fibration from a surface onto a curve is said to be of 
 {\em genus $1$}
 if its generic fiber is an integral curve of arithmetic genus $1$.
 In case the generic fiber is smooth,
 the fibration is called {\em elliptic}, and {\em quasi-elliptic} otherwise.
 Moreover, if the fibration admits a section, it is called 
 {\em Jacobian}, and
 a choice of section, referred to as the {\em zero section}, is part of the data.
\end{Definition}

In characteristic different from $2$ and $3$, the generic fiber of a genus $1$
fibration is automatically smooth by \cite{bm3}, and thus, an elliptic fibration.

\subsection{Families of torsors arising from formal Brauer groups}
\label{subsec:moving torsor general}
For future applications, we extend our setup in this subsection
and work with Jacobian genus $1$ fibrations from surfaces that 
are not necessarily K3.
We follow the setup of the articles \cite{artin swinnerton}
and \cite{artin} by Artin and Swinnerton-Dyer.
Let
$$f\,:\,X\,\to\, Y$$ 
be a relatively minimal 
(that is, there are no $(-1)$-curves in the fibers)
Jacobian genus $1$ fibration, where 
$X$ is a surface, and $Y$ is a curve, both smooth and proper
over an algebraically closed field $k$.
Contracting those $(-2)$-curves in the fibers of $f$ that do not
intersect the zero section, we obtain the {\em Weierstra\ss\ model}
$$
  f'\,:\,X'\,\to\, Y.
$$
If $f$ has reducible fibers, then $X'$ has rational double point singularities.
We denote by $A\subseteq X'$ the smooth locus of $X'$.
As explained in \cite[Section 1]{artin swinnerton},
$A$ has a unique structure $\oplus$ of group scheme over $Y$:
namely, if $P_1, P_2$ are sections of $A$ over $Y$, then they are Cartier divisors,
and $P_1\oplus P_2$ is the zero locus of a non-zero section of
$\OO_{X'}(P_1+P_2-Z)$, where $Z$ denotes the zero section.
In case $f$ is an elliptic fibration, we have the following interpretation
in terms of N\'eron models: the smooth locus 
of $X$ over $Y$ is the N\'eron model
of its generic fiber, and $A$ is its identity component.

Next, let $S$ be the formal spectrum of a local, 
Noetherian, and complete $k$-algebra with residue field $k$. 
We want to classify families of torsors under $A$, parametrized by $S$,
such that the special fiber is the trivial $A$-torsor.
That is, we consider Cartesian diagrams of algebraic spaces
\begin{equation}
\label{eq:torsor diagram}
  \xymatrixcolsep{3pc}
  \xymatrix{
      A \ar[r] \ar[d] & \oA \ar[d] \\
      Y \ar[r] \ar[d] & Y\times_k S \ar[d] \\
       \Spec k \ar[r] & S  }
\end{equation}
In order to classify such {\em moving torsors}, we recall that 
Artin and Mazur \cite{artin mazur} studied
the functors on local Artinian $k$-algebras with residue field $k$
$$
 \begin{array}{ccccc}
  \Phi_{X/k}^i &:& (\mbox{Art}/k) & \to & (\mbox{Abelian groups})\\
     & & R &\mapsto&\ker\left( \Het{i}(X\times_k\Spec R, \GG_m)\,\to\, \Het{i}(X,\GG_m) \right)
 \end{array}
 $$
see also Section \ref{subsec: brauer}.
We now furthermore assume that $\Phi_{X/k}^2$ is pro-representable
by a formal group law, which is then called the {\em formal Brauer group} and denoted
$\widehat{{\rm Br}}(X)$,
Next, let us recall that there exists a 
Grothendieck--Leray spectral sequence 
$$
   E_2^{i,j}\,:=\,\Het{i}(Y, \,R^jf'_\ast\,\GG_m) \,\Longrightarrow\,
    \Het{i+j}(X',\, \GG_m)\,.
$$
As Artin explained in \cite[Section 2]{artin}, the formal structure of 
$\Het{2}(X,\GG_m)$ is that of $\Het{1}(Y,\Pic_{X'/Y})$.
Using the zero section of $f'$, we identify $\Pic_{X'/Y}^0$ with $A$, 
and then, it is not difficult to see that moving torsors 
are closely related to the formal Brauer group.
More precisely, we have the following result.

\begin{Proposition}
   \label{prop:torsor families}
   We keep the notations and assumptions.
   Let $S:=\Spf R$, where $(R,\idealm_R)$ is a local, Noetherian, and 
   complete $k$-algebra with residue field $k$.
   Let $n\geq1$ be an integer.
   \begin{enumerate}
      \item Formal families of $A$-torsors $\oA\to Y\times_k S$,
        whose special fiber is the trivial $A$-torsor, are classified
        by the $R$-valued points
        $$
           \widehat{{\rm Br}}(X)(R)
        $$
        of the formal Brauer group of $X$.
      \item The compactification $A\subseteq X'$ extends to 
        a compactification $\oA\subseteq \oX'$,
        and the formal family $\oX'\to Y\times_k S$ is algebraizable.
      \item Moreover, $n$-torsion elements of $\widehat{{\rm Br}}(X)(R)$
        correspond to families as in (1)
        such that there exists 
        a degree $n$ section of $\Pic_{\oX'/Y\times_k S}$ over $Y\times_k S$.
  \end{enumerate}
\end{Proposition}

\prf
First, we use the zero section of $f'$ to identify
$\Pic^0_{X'/Y}$ with $A$.
Then, as explained at the beginning of \cite[Section 2]{artin}
and in \cite[Proposition (2.1)]{artin},
the formal structures of 
$\Het{2}(X,\GG_m)$ and $\Het{2}(X',\GG_m)$ 
are that of  $\Het{1}(Y,\Pic_{X'/Y})$.
That is, by definition of $\Phi_{X/k}^2$ and its pro-representability 
assumption, we have
$$
\widehat{{\rm Br}}(X)(R)\,=\,
\ker\left(
\Het{1}(Y\times_k S, A)\,\stackrel{{\rm res}}{\longrightarrow}\,\Het{1}(Y,A)
\right)\,,
$$
where ${\rm res}$ denotes restriction.
But then, elements of the right hand side classify $A$-torsors over
$Y\times_k S$, whose restriction to the special fiber is trivial.
This shows claim (1).

Next, we show compactification of $\oA$.
We set $R_m:=R/\idealm^m$ and $S_m:=\Spec R_m$.
To simplify notations, we denote by $-_{S_m}$ the trivial
product family $-\times_{\Spec k}S_m$.
By induction on $m$, we may assume that we have already extended
the compactification $A\subseteq X'$ to some $\oA_{S_m}\subseteq\oX'_m$.
Blowing up the boundary, 
we obtain a compactification $\oA_{S_m}\subseteq\oY_m$,
whose boundary is a Cartier divisor.
As explained in \cite[Section 2.1]{clo}, this latter
compactification can be extended to a compactification 
$\oA_{S_{m+1}}\subseteq \oY_{m+1}$.
Blowing down $\oY_{m+1}$ to 
$\oX'_m$ 
(see, for example, in \cite[Theorem 3.1]{cynk van straten}),
we obtain a compactification $\oA_{S_{m+1}}\subseteq\oX'_{m+1}$,
which extends $\oA_{S_m}\subseteq\oX'_m$.
Passing to the limit, we obtain the desired compactification
$\oA\subseteq\oX'$.

Multiplication by $n$ induces a morphism $A\to A$ of group schemes
over $Y$, and thus, a morphism $\tau_n:\Het{1}(Y,A)\to\Het{1}(Y,A)$.
From the discussion at the end of \cite[Section 1]{artin swinnerton}
it follows that an element in the kernel of
$\tau_n$ corresponds to an $A$-torsor over $Y$ such that there exists
a section of $\Pic_{X'/Y}$ over $Y$ of degree $n$.
The same holds true with $Y$ replaced by $Y_S$,
and thus, $n$-torsion elements of $\widehat{{\rm Br}}(X)(R)$
correspond to formal families of $A$-torsors over $Y_{S}$
that become trivial over the special fiber,
such that there exists a degree-$n$ section of 
$\Pic_{\oX'/Y_{S}}$  over  $Y_{S}$.
This shows claim (3).

It remains to show algebraization.
By the established assertion (3), there exists
a degree-$n$ section $\overline{\cal L}$ of  $\Pic_{\oX'/Y_S}$ over $Y_S$.
Since $Y$ is a curve over an algebraically closed field, we have ${\rm Br}(Y)=0$ 
by Tsen's theorem.
Since $H^2(Y,\OO_Y)=0$, we have $\widehat{{\rm Br}}(Y)=0$, which implies
${\rm Br}(Y_S)=0$,
and we obtain a short exact sequence
$$
0\,\to\,\Pic(Y_S)\,\to\,\Pic(\oX')\,\to\,
H^0(Y_S,\,\Pic_{\oX'/Y_S})\,\to\,
\underbrace{{\rm Br}(Y_S)}_{=0}\,\to\,...
$$
In particular, $\overline{\cal L}$ 
lifts to some $\widetilde{\cal L}\in\Pic(\oX')$.
Next, let $E\in\Pic(\oX')$ be the class of a fiber, and then, for every integer $m$,
we define ${\cal M}_m:=\widetilde{{\cal L}}\otimes\OO_{\oX'}(mE)$.
Since every integral curve on $X'$ is either a fiber or a multisection of the fibration,
it follows that the restriction ${\cal M}_m|_{X'}$ has positive
intersection with every integral curve on $X'$ if $m\gg0$.
Moreover, for $m\gg0$, the self-intersection of ${\cal M}_m|_{X'}$ 
is positive.
Thus, by the Nakai--Moishezon criterion, 
for $m\gg0$, the restriction of ${\cal M}_m$ to $X'$ is an ample 
invertible sheaf.
Therefore,  the formal family $\oX'$ is algebraizable 
by Grothendieck's existence theorem, which establishes claim (2).
\qed\medskip

Before proceeding, let us recall a couple of facts about
commutative formal group laws, and refer,
for example, to 
\cite{zink} for details:
if $\widehat{F}$ is a commutative formal group law of dimension $d$
over a field of characteristic zero, then there exists a unique 
strict isomorphism to $\widehat{\GG}_a^d$, the logarithm of $\widehat{F}$. 
On the other hand, if $\widehat{F}$ is defined over
a field of positive characteristic $p$, then
there exists a short exact sequence of commutative formal group laws
$$
 0\,\to\,\widehat{F}^{\rm u}\,\to\,\widehat{F}\,\to\,\widehat{F}^{\rm bt}\,\to\,0,
$$
where $\widehat{F}^{\rm u}$ is unipotent and $\widehat{F}^{\rm bt}$ is $p$-divisible
\cite[Theorem 5.36]{zink}.
We recall that a formal group law $\widehat{F}$ is called {\em $p$-divisible}
if multiplication by $p$ is an isogeny, and then, there exists
an integer $m\geq1$ such that the $m$-fold Frobenius 
$$
   {\rm Fr}^{m}\,:\,\widehat{F}\,\to\,\widehat{F}^{(p^m)}
$$
factors through multiplication by $p$.
The minimal $m$, for which such a factorization exists, is called
the {\em height} of $\widehat{F}$.
On the other extreme, multiplication by $p$ on $\widehat{\GG}_a$
is zero and thus, this formal group law is of infinite height.
More generally, if $\widehat{F}$ is {\em unipotent}, then there exists
an increasing sequence of formal subgroup laws 
$0=\widehat{F}_0\subset...\subset\widehat{F}_r=\widehat{F}$
such that successive quotients are isomorphic to $\widehat{\GG}_a$,
see \cite[Theorem 5.37]{zink}.

This recalled, we have the following statement about formal group laws only,
which we need to ensure the existence of non-trivial moving torsor families
over $\Spec k[[t]]$.

\medskip\begin{Lemma}
  \label{lemma:formal groups}
  Let $\widehat{F}$ be a formal group law over $k$, and let 
  $(R,\idealm_R)$ be a local, Noetherian, and complete $k$-algebra with residue
  field $k$.
  \begin{enumerate}     
    \item If $p$ does not divide $n$, or $R$ is reduced and $\widehat{F}$ 
      is a $p$-divisible formal group law, then 
      $$
         \widehat{F}(R)[n] \,=\,0.
       $$
    \item If $R$ is reduced and $\idealm_R\neq0$, then
      $$
         \widehat{F}(R)[p]\,\neq\,0\mbox{ \quad }\Leftrightarrow\mbox{ \quad }
         \widehat{F}^{\rm u}\,\neq\,0\,,
      $$
      that is, if and only if $\widehat{F}$ is not $p$-divisible.
   \end{enumerate}
\end{Lemma}

\prf
If $p\nmid n$, then multiplication by $n$ is injective, and thus, $\widehat{F}(R)[n]=0$.
If $\widehat{F}$ is $p$-divisible, say, of finite height $h$, then the $h$-fold
Frobenius factors through multiplication by $p$.
Since Frobenius is injective on $R$-valued points of $\widehat{F}$ for
reduced $R$, this implies $\widehat{F}(R)[p]=0$
and establishes claim (1).

If $\widehat{F}^{\rm u}=0$ and $R$ is reduced, then 
$\widehat{F}(R)[p]=0$ by assertion (1).
Conversely, if $\widehat{F}$ is not $p$-divisible, then 
$\widehat{\GG}_a\subseteq\widehat{F}$.
Since $\widehat{\GG}_a(R)[p]=\idealm_R\neq0$, we find
$\widehat{F}(R)[p]\neq0$, which establishes claim (2).
\qed\medskip

Proposition \ref{prop:torsor families} and Lemma \ref{lemma:formal groups} 
give us a criterion for the existence of
families of $A$-torsors, and ensure compactification and algebraization.
The following result gives desingularization, and compares 
the formal Brauer groups and the Picard groups
of special and geometric generic fiber -- we note that the latter implies that the family
has non-trivial moduli.

\begin{Proposition}
  \label{prop:picard index p}
  We keep the notations and assumptions of Proposition \ref{prop:torsor families}
  and assume that $\widehat{{\rm Br}}(X)$ is not $p$-divisible.
  Set $R:=k[[t]]$ and $S:=\Spec R$, and let
  $$
     \oA \,\to\, Y\times_{\Spec k} S \,\to\, S
  $$ 
  be a family of $A$-torsors associated to a non-trivial $p$-torsion 
  element of $\widehat{{\rm Br}}(X)(R)$
  as in Proposition \ref{prop:torsor families}.
  \begin{enumerate}
    \item Let $\oX'\to S$ be a compactification as asserted
      in Proposition \ref{prop:torsor families}.
      Then, after possibly replacing $S$ by a finite flat cover,
      there exists a simultaneous resolution of the singularities
      $$
      \xymatrixcolsep{3pc}
      \xymatrix{
         X \ar[r] \ar[d] & \oX \ar[d] \\
         Y \ar[r] \ar[d] & Y\times_{\Spec k} S \ar[d] \\
         \Spec k \ar[r] & S  }
      $$
      which is smooth and projective over $S$.
      \item Let $\overline{\eta}\in S$ be the geometric generic point.
      Then, there exists an isomorphism of formal group laws
      $$
          \widehat{{\rm Br}}(\oX_{\overline{\eta}})\,\iso\,
          \widehat{{\rm Br}}(X) \,\otimes_k\, \kappa(\overline{\eta}).
      $$
      In particular, if $X$ is a supersingular K3 surface, then so
      is $\oX_{\overline{\eta}}$.      
    \item 
    Specialization induces a short exact sequence of Picard groups
    $$
          0\,\to\,\Pic(\oX_{\overline{\eta}})
          \,\to\,\Pic(X)\,\to\,\ZZ/p\ZZ\,\to\,0\,,
     $$
      whose cokernel is generated by the class of the
      zero section of $f:X\to Y$.
      \end{enumerate}
\end{Proposition}

\prf
First, $\oX'\to S$ is a flat family of surfaces, 
whose special fiber $X'$ has at worst rational double points as singularities.
Thus, also the generic fiber has at worst
rational double points as singularities
by \cite[Proposition 6.1]{liedtke noether},
and thus, after possibly base-changing to a finite flat extension of $S$,
there exists a simultaneous resolution of singularities $\oX\to S$
by the main result of \cite{artin resolution}.
This establishes claim (1) except for the projectivity statement.

Claim (2) follows from \cite[Proposition (2.1)]{artin}.

To establish claim (3), we note that there exists a commutative diagram with 
exact rows, whose vertical arrows are restriction maps:
$$
\begin{array}{ccccccccc}
 0 & \to& \Pic(Y_{\overline{\eta}}) &\to& \Pic (\oX'_{\overline{\eta}}) &\to& 
    H^0(Y,\, \Pic_{\oX'_{\overline{\eta}}/Y_{\overline{\eta}}}) &\to& 0 \\
  &&\uparrow&&\uparrow&&\uparrow\\
 0 & \to& \Pic(Y_S) &\to& \Pic ( \oX' ) &\to& 
    H^0(Y_S,\, \Pic_{\oX'/Y_S}) &\to& 0\\
   &&\downarrow&&\downarrow&&\downarrow\\
 0 & \to& \Pic(Y) &\to& \Pic (X') &\to& 
    H^0(Y,\, \Pic_{X'/Y}) &\to& 0
\end{array}
$$
see \cite[(2.2)]{artin swinnerton}, or the proof of 
Proposition \ref{prop:torsor families}.
Replacing $S$ by a finite flat cover, we may and will assume
that the Picard groups of $\oX'_\eta$ and $\oX'_{\overline{\eta}}$ are isomorphic.
Next, it follows from \cite[Proposition (1.6)]{artin swinnerton} that there exists
a commutative diagram of group algebraic spaces
over $Y$, $Y_S$ and $Y_{\overline{\eta}}$, respectively:
$$
\begin{array}{ccccccccc}
0 &\to& A_{\overline{\eta}} &\to& \Pic_{\oX_{\overline{\eta}}'/Y_{\overline{\eta}}} 
    &\to& \underline{\ZZ}_{Y_{\overline{\eta}}} &\to& 0\\
   &&\uparrow&&\uparrow&&\uparrow\\   
  0 &\to& A_S &\to& \Pic_{\oX'/Y_S} &\to& 
    \underline{\ZZ}_{Y_S} &\to& 0\\    
 &&\downarrow&&\downarrow&&\downarrow\\
  0 &\to& A &\to& \Pic_{X'/Y} &\to& 
    \underline{\ZZ}_Y &\to& 0
\end{array}
$$
The class of the zero section $Z$ of $X'\to Y$ in $\Pic(X')$ 
defines a splitting of the bottom row.
By Proposition \ref{prop:torsor families}, there exists a 
degree-$p$ section of $\Pic_{\oX'/Y_S}$ over $Y_S$.
Thus, taking global sections in the previous diagram, we conclude that
the image of $H^0(Y_S,\Pic_{\oX'/Y_S})$ inside 
$H^0(Y_S,\underline{\ZZ}_{Y_S})$ is of index $1$ or $p$.
However, this index cannot be equal to $1$, since $\oA\to Y_S$
is a non-trivial family of $A$-torsors.
Combining these observations and the two commutative diagrams,
we arrive at a short exact sequence of Abelian groups
$$
0\,\to\,\Pic(\oX'_{\overline{\eta}})\,\to\,\Pic(X')\,\to\,\ZZ/p\ZZ\,\to\,0,
$$
where the cokernel is generated by the class of $Z$.

Next, $\Pic(X)$ is generated by the exceptional divisors of
the contraction morphism $\nu:X\to X'$ and 
$\nu^\ast\Pic(X')$, and we have a similar statement 
for $\Pic(\oX_{\overline{\eta}})$.
Since $\oA\to Y_S$ is a family of $A$-torsors, 
and the special fiber $A$ has no multiple fibers, neither has the generic
fiber, and thus, the singular fibers do not change their type by
\cite[Theorem 5.3.1]{cd}.
In particular, $\oX'_{\overline{\eta}}$ and $X'$ 
have the same types of rational double points.
From this, we deduce that also the cokernel of the specialization
homomorphism $\Pic(\oX_{\overline{\eta}})\to\Pic(X)$ 
is cyclic of order $p$, generated by the class of $Z$, 
which establishes claim (2).

Finally, if ${\cal L}$ is an ample invertible sheaf on $X$, then
${\cal L}^{\otimes p}$ extends to $\oX$, which shows that 
$\oX\to S$ is projective.
\qed\medskip

To understand the geometry of moving $A$-torsors better, and to construct
purely inseparable multisections, we now inspect the generic fiber 
more closely.
Let $\xi\in Y$ be the generic point, and restrict a family
$\oA\to Y\times_k S$ as in Proposition \ref{prop:torsor families}
to $\xi\times_k S$, that is, we consider
$$
   \oA_\xi \,:=\, \oA\times_{(Y\times_k S)}(\xi\times_k S) 
   \,\to\, \xi\times_k S\,.
$$
This is a family of torsors under $A_\xi:=A\times_Y\xi$ over $\xi\times_k S$.
Let us recall that the relative Frobenius $F:A\to A^{(p)}$ is a
morphism of group schemes over $Y$, whose kernel 
$A[F]$ is a finite, flat, and infinitesimal group scheme of length $p$ over $Y$.

\begin{Proposition}
  \label{prop:inseparable multisection}
  We keep the notations and assumptions of Proposition \ref{prop:torsor families}
  and assume that $\widehat{{\rm Br}}(X)$ is not $p$-divisible.
  Set $R:=k[[t]]$ and $S:=\Spec R$, and let
  $$
     \oA \,\to\, Y\times_{k} S \,\to\, S
  $$ 
  be a family of $A$-torsors
  associated to a non-trivial $p$-torsion element of $\widehat{{\rm Br}}(X)(R)$
  as in Proposition \ref{prop:torsor families}.
  Then, after possibly replacing $S$ by some finite flat cover
  \begin{enumerate}     
    \item There exists a degree $p$ multisection 
      $D_\xi\subset\oA_\xi$
      such that the induced morphism
      $$
         D_\xi\,\to\, \xi\times_{k} S
      $$
      is finite, flat and radiciel of degree $p$.
    \item More precisely, $D_\xi\to\xi\times_k S$
      is a family of $A_\xi[F]$-torsors.
      Thus, we obtain an isomorphism
       $$
          \oA_\xi \,\iso\,
          \left(A_\xi \times_{(\xi\times_{k} S)} D_\xi\right) \,/\, A_\xi[F]\,,
       $$
       that is, a description of this family of $A_\xi$-torsors as $A_\xi[F]$-twist.
  \end{enumerate}
\end{Proposition}

\prf
If $Z$ denotes the zero section of $f':X'\to Y$, 
then $\OO_{X'}(pZ)$ extends to some invertible sheaf $\cal M$
on $\oX'$ by Proposition \ref{prop:picard index p}.
Since $\cal M$ has degree $p$ on each fiber,
Riemann--Roch implies that 
it has vanishing first cohomology and a $p$-dimensional space
of global sections on each fiber.
Thus, $\overline{\cal Q}:=(\overline{f}'_\ast{\cal L})^\vee$ 
is a locally free $\OO_{Y\times_k S}$-module of rank $p$.
This introduced, we recall (for example, from \cite[Section 3]{kleiman picard})
that relative effective Cartier divisors $D\to Y_S$ such
that $\OO_{X'}(D)$ is isomorphic to ${\cal M}$ modulo invertible
sheaves coming from $Y_S$, correspond to sections of
$\PP(\overline{\cal Q})\to Y_S$.
Since the fiber of $\PP(\overline{\cal Q})\to Y_S$ over $\xi\times_k S$
is isomorphic to $\PP^{p-1}_{\xi\times_k S}$, this already shows the existence of 
degree-$p$ multisections of $\oA_\xi\to \xi\times_kS$.

Next, we show that $A[F]$ acts on $\PP(\overline{\cal Q})\to Y_S$
by translation of relative effective Cartier divisors:
the $A$-action on $\oA$ induces an $A$-action on the symmetric
product ${\rm Sym}^p(\oA)$, which we identify 
with the space of relative effective Cartier divisors of $\oX'\to Y_S$ 
of degree $p$ (see, for example, \cite[Remark 9.3.9]{kleiman picard}).
Inside this latter space we have $\PP(\overline{\cal Q})\to Y_S$, which is
stable under the induced $A[F]$-action.


We now determine the schematic fixed point locus of the $A_\xi[F]$-action on 
$\PP^{p-1}_{\xi\times S}$.
On geometric fibers, a fixed point is of the form $pP$.
Thus, depending on the $p$-torsion subgroup scheme of $A_\xi$, 
the fixed point locus is either Artinian of length $p^2$ (if $f'$ is elliptic), or
it is a curve (if $f'$ is quasi-elliptic).
In any case, the fixed point locus is flat over $\xi\times S$.
Thus, after possibly replacing $S$ by a finite flat cover, there exists
an $A_\xi[F]$-invariant section $D_\xi$ of $\PP^{p-1}_{\xi\times_k S}\to\xi\times_k S$ that specializes to $pZ$.
By construction, $D_\xi\to \xi\times_k S$ is a family of $A_\xi[F]$-torsors, and in particular,
finite, flat, and radiciel of degree $p$ over $\xi\times_k S$, 
which establishes claim (1).
Since the base-change of $\oA_\xi$ to $D_\xi$ trivializes the torsor,
we obtain a description of $\oA_\xi$ as $A_\xi[F]$-twist, which
establishes claim (2).
\qed\medskip

We now summarize the results on moving $A$-torsors obtained so far and 
use a purely inseparable degree-$p$ multisection as established in the previous 
proposition to show that special and generic fiber of a family 
of moving $A$-torsors are related by a purely inseparable isogeny.
We note that the following theorem is the technical heart of this article.

\begin{Theorem}
   \label{thm:moving torsors general}
  We keep the notations and assumptions of Proposition \ref{prop:torsor families}
  and assume that $\widehat{{\rm Br}}(X)$ is not $p$-divisible.
  Let $R:=k[[t]]$ and $S:=\Spec R$ and let
  $$
     \oA \,\to\, Y\times_{\Spec k} S \,\to\, S
  $$ 
  be a family of $A$-torsors
  associated to a non-zero $p$-torsion element of $\widehat{{\rm Br}}(X)(R)$
  as in Proposition \ref{prop:torsor families}.
  Then, after possibly replacing $S$ by a finite and flat cover,
  \begin{enumerate}     
    \item There exists a compactification and desingularization
      of $\oA\to Y\times_k S$ to a 
      smooth and projective family 
      $$
         \oX\,\to\, Y\times_{\Spec k} S\,\to\,S
      $$
      with special fiber $X$.
    \item Specialization induces a short exact sequence
      $$
         0\,\to\,\Pic(\oX_{\overline{\eta}})\,\to\,\Pic(X)\,\to\,\ZZ/p\ZZ\,\to\,0,
      $$ 
      and we have an isomorphism
      $$
         \widehat{{\rm Br}}(\oX_{\overline{\eta}})
         \,\iso\, 
         \widehat{{\rm Br}}(X)\otimes_k \kappa(\overline{\eta}),
      $$
      where $\overline{\eta}$ denotes the geometric generic point of $S$.
     \item There exist a morphism and a rational map
      $$
       \xymatrixcolsep{3pc}
       \xymatrix{
        & Y \ar[dl] \ar@{.>}[dr] \\
        X\times_{\Spec k}\eta & & \oX_\eta
        }
     $$
     both of which are generically finite and purely inseparable of
     degree $p$.   
    \item  
      There exist rational maps
     $$
         \oX_{\overline{\eta}} \,\dashrightarrow\,X\times_{\Spec k} \overline{\eta} 
         \,\dashrightarrow\, \oX_{\overline{\eta}},
     $$
     both of which are generically finite and purely inseparable
     of degree $p^2$.
     Thus, $\oX_{\overline{\eta}}$ and $X\times_k\overline{\eta}$ are purely
     inseparably isogenous of height $2$.
    \end{enumerate}
\end{Theorem}

\prf
We established claims (1) and (2) in Proposition \ref{prop:torsor families} and 
Proposition \ref{prop:picard index p}.

After possibly replacing $S$ by a finite flat cover, there exists
a purely inseparable degree-$p$
multisection $D_\xi\subset \oA_\xi$  by
Proposition \ref{prop:inseparable multisection}, and we denote by
$D$ its closure in $\oX_\eta$.
Since $D_\xi\to\xi\times_k S$ is finite, flat and radiciel of degree $p$, 
the same is true for $D\to Y\times_k \eta$.
Base changing to $D\to Y\times_k \eta$ trivializes
the compactified family of $A$-torsors generically,
and therefore, we obtain a diagram
$$
   X\times_{\Spec k} \eta  \,\longleftarrow\,
   (X\times_{\Spec k} \eta)\times_{(Y\times_k \eta)} D
   \,\stackrel{\iso}{\dashrightarrow}\,
   \oX_\eta\times_{(Y\times_k\eta)}D\,\longrightarrow\,\oX_\eta\,,
$$
where the morphisms on the left and right are purely inseparable of degree $p$,
and the rational map in the middle is birational.
This establishes claim (3).

Let $F_{\overline{\eta}}:\oX_{\overline{\eta}}\to\oX_{\overline{\eta}}^{(p)}$
be the $\overline{\eta}$-linear Frobenius morphism,
choose an $\overline{\eta}$-linear isomorphism 
$\oX_{\overline{\eta}}^{(p)}\iso\oX_{\overline{\eta}}$, and note that
$F_{\overline{\eta}}$  factors through $\oX_\eta\times_{(Y\times_k\eta)}D$,
see also Section \ref{subsec: isogeny}.
From this, we obtain a composition
$$
\oX_{\overline{\eta}}\,\to\,
\oX_{\overline{\eta}}\times_{(Y\times_k\overline{\eta})}D\,\dashrightarrow\,
X\times_{\Spec k}\overline{\eta},
$$
which is a rational map of varieties over $\overline{\eta}$, which is generically finite and 
purely inseparable of degree $p^2$, that is, a purely inseparable isogeny of height $2$.
This establishes claim (4).
\qed\medskip

\subsection{Families of supersingular K3 surfaces}
In this subsection we specialize to K3 surfaces.
We recall from Section \ref{subsec: brauer} that 
the formal Brauer group $\widehat{{\rm Br}}(X)$
of a K3 surface $X$ is a one-dimensional formal group law.
In particular, we have the following equivalences:
$$
\widehat{{\rm Br}}(X)^{\rm u}\,\neq\,0
\,\Leftrightarrow\, 
h\left(\widehat{{\rm Br}}(X)\right)\,=\,\infty
\,\Leftrightarrow\, \mbox{ $X$ is supersingular} .
$$
Thus, by Proposition \ref{prop:torsor families} and Lemma \ref{lemma:formal groups},
non-trivial families of moving torsors
over $\Spec k[[t]]$ associated to a Jacobian (quasi-)elliptic K3 surface 
can exist only for supersingular K3 surfaces,
which renders precise Artin's remark: ``The unusual phenomenon of continuous families of homogeneous spaces 
occurs only for supersingular surfaces'' \cite[footnote (2) on p. 552]{artin}.
 The next proposition rephrases Theorem \ref{thm:moving torsors general}
in terms of supersingular K3 surfaces.

\begin{Proposition}
 \label{prop: deformation}
 Let $X\to\PP^1$ be a Jacobian (quasi-)elliptic fibration on a supersingular K3 surface over $k$.
 Then, there exists a smooth and projective family of supersingular elliptic K3 surfaces
 with non-trivial moduli
 $$
   \oX\,\to\,\PP^1_S\,\to\,S,\mbox{ \quad where \quad } S\,:=\,\Spec k[[t]],
 $$
 whose special fiber is $X\to\PP^1$ and
 that has the following properties:
 \begin{enumerate}
  \item The Artin invariant of the geometric generic fiber satisfies
  $$
    \sigma_0( \oX_{\overline{\eta}}) \,=\, \sigma_0(X)\,+\,1.
  $$ 
   \item There exist purely inseparable isogenies of height $2$,
   that is, dominant, rational, and generically finite maps
     $$
         \oX_{\overline{\eta}} \,\dashrightarrow\,X\times_{\Spec k} \overline{\eta} 
         \,\dashrightarrow\, \oX_{\overline{\eta}},
     $$
   whose composition is twice the $\overline{\eta}$-linear Frobenius morphism.
 \end{enumerate}
\end{Proposition}

\prf
By Theorem \ref{thm:moving torsors general},
the index of $\Pic(\oX_{\overline{\eta}})$ in
$\Pic(X)$ is equal to $p$, and thus,
claim (1) follows from the definition of the Artin invariant.
In particular, since the Artin invariants of $X$ and $\oX_{\overline{\eta}}$ 
differ, the family has non-trivial moduli.
The remaining assertions are  explicitly stated in 
Theorem \ref{thm:moving torsors general}.
\qed\medskip

In characteristic $p\geq5$, 
supersingular K3 surfaces do not degenerate, that is,
have potential good reduction, 
by a theorem of Rudakov and Shafarevich \cite{rudakov shafarevich degeneration}.
Thus, the family over $\Spec k[[t]]$ described in the previous
proposition can be spread out to a {\em smooth} family of supersingular K3 surfaces 
over a smooth and {\em proper} curve.
More precisely, we have the following result.

\begin{Proposition}
 \label{prop: spread out} 
 If $p\geq5$ and under the assumptions of Proposition \ref{prop: deformation},
 there exist a smooth projective curve $C$ over $k$,
 a closed point $0\in C$, and
 a smooth projective family of supersingular K3 surfaces 
 $$
    \oY \,\to\,C
 $$
 with the following properties:
 \begin{enumerate}
   \item After possibly replacing $S$ by a finite flat cover, 
      $\oX\to S$ is the fiber over the completed local ring
      $\widehat{\OO}_{C,0}$. In particular, $X$ is the fiber over $0$.
  \item Specialization induces an embedding 
      $$
         \Pic(\oY_{\overline{\eta}})\,\subset\,\Pic(X),
      $$
      which is of index $p$.
      More precisely, if $E$ denotes a fiber of $X\to\PP^1$,
      and $Z$ the zero-section, then the classes of
      $E$ and $pZ$ extend to $\Pic(\oY_{\overline{\eta}})$.
  \item Let $c\in C$ be a point such  
      that the geometric fiber $\oY_{\overline{c}}$
      has Artin-invariant $\sigma_0(X)+1$.
      Then, specialization of $E$ to $\oY_{\overline{c}}$
      gives rise to a non-Jacobian elliptic fibration.
      Moreover, there exists a purely inseparable degree-$p$ multisection 
      $D_{\overline{c}}$ on $\oY_{\overline{c}}$, which of class 
      $pZ+kE$ for some $k\geq2$.
   \item Under the assumptions of (3), there exist purely inseparable
      isogenies
      $$
         \oY_{\overline{c}}\,\dashrightarrow\, X\times_{\Spec k} \Spec \kappa(\overline{c})
         \,\dashrightarrow\,\oY_{\overline{c}}\,,
       $$
      both of which are of height $2$.
 \end{enumerate}
\end{Proposition}

\prf
By Artin's approximation theorem \cite[Theorem 1.6]{artin approx},
the family $\oX\to S$ can be defined over a 
$k$-algebra of finite type.
From there, we spread it out 
to a projective family $\oY\to C$, where $C$ is
a smooth projective curve over $k$.
We denote by $0\in C$  the point such that the family over the
completed ring $\widehat{\OO}_{C,0}$ is $\oX$.
Since supersingular K3 surfaces in characteristic $p\geq5$ 
have potential good reduction by \cite{rudakov shafarevich degeneration}, 
we may assume, after possibly replacing $C$ by a finite flat
cover, that $\oY\to C$ is a smooth projective
family of supersingular K3 surfaces.
This establishes claim (1).

We have a family of elliptic fibrations $\oX\to\PP^1\times_k S\to S$
(since $p\geq5$,  the fibrations cannot be quasi-elliptic).
In particular, the class of $E$ extends from $X$ to $\oY_{\overline{\eta}}$,
which, together with Proposition \ref{prop:picard index p} establishes claim (2).

Now, let $c\in C$ be a point such that 
$\sigma_0(\oY_{\overline{c}})=\sigma_0(\oY_{\overline{\eta}})$.
Then, specialization induces an isomorphism 
$\Pic(\oY_{\overline{\eta}})\iso\Pic(\oY_{\overline{c}})$.
The elliptic fibration on $\oY_{\overline{\eta}}$ specializes to an elliptic
fibration on $\oY_{\overline{c}}$.
However, this latter fibration cannot be Jacobian,
for otherwise there would exist a section, whose class 
would extend to $\oY_{\overline{\eta}}$, and which would
give rise to a section of the original elliptic fibration of $\oY_{\overline{\eta}}$, 
a contradiction.
Since $\oX\to\PP^1\times_k S$ is a family of $A$-torsors,
also the Jacobian fibration associated to $\oY_{\overline{c}}\to\PP^1$
is $X\to\PP^1$.

Next, the degree-$p$ multisection $D\subset\oY_{\overline{\eta}}$ 
from the proof of Theorem \ref{thm:moving torsors general}
specializes to a degree-$p$ multisection $D_{\overline{c}}\subset\oY_{\overline{c}}$.
Now, $D_{\overline{c}}$ must be an integral curve, for otherwise,
a linear combination of $D_{\overline{c}}$ and $(D_{\overline{c}})_{\rm red}$
would give rise to a relative invertible sheaf on $\oY_{\overline{c}}\to \PP^1$ of degree $1$,
contradicting the fact that this fibration is not Jacobian.
Since the class of $D$ on $\oY_{\overline{\eta}}$ is equal to 
$pZ$ modulo fiber classes, it must be of class $pZ+kE$ for some integer $k$,
and similarly for $D_{\overline{c}}$.
Since integral curves on K3 surfaces have self-intersection number at least $-2$,
we compute $k\geq2$.
And finally, since $D$ is purely inseparable of degree $p$ 
over the base, the same is true for its specialization $D_{\overline{c}}$.
This establishes claim (3).

Having a non-Jacobian elliptic fibration $\oY_c\to\PP^1$ with a purely
inseparable degree-$p$ multisection $D_{\overline{c}}$, whose associated Jacobian 
fibration is $X\to\PP^1$, the same arguments for the proof of assertion (4) of Theorem \ref{thm:moving torsors general} 
also show that $\oY_{\overline{c}}$ is related to $X$ by a purely inseparable isogeny 
of height $2$.
This establishes claim (4).
\qed\medskip

\subsection{Jacobian elliptic fibrations on supersingular K3 surfaces}
In order to use Proposition \ref{prop: deformation}, we have to show the
existence of Jacobian elliptic fibrations on supersingular K3 surfaces.
For example, a supersingular K3 surface with Artin invariant $\sigma_0=10$
cannot possess such a fibration, for otherwise Proposition \ref{prop: deformation}
would produce a supersingular K3 surface with $\sigma_0=11$,
which is impossible.
The next proposition shows that this is the only restriction.

\begin{Proposition}
  \label{prop:elliptic fibrations}
  Let $X$ be a supersingular K3 surface with Artin invariant $\sigma_0$
  in characteristic $p\geq5$. 
  \begin{enumerate}
    \item If $\sigma_0\leq9$, then $X$ admits a Jacobian elliptic fibration.
    \item If $\sigma_0=10$, then $X$ does not admit a Jacobian elliptic fibration.
  \end{enumerate}
\end{Proposition}

\begin{Remark}
  Assertion (2) was already shown by 
  Ekedahl and van~der~Geer \cite[Proposition 12.1]{Ekedahl van der Geer},
  as well as by Kond\={o} and Shimada 
  \cite[Corollary 1.6]{kondo shimada}, but using different methods.
\end{Remark}

\prf
We have shown claim (2) in the lines before this proposition.

By \cite[Section 1]{rudakov shafarevich}, the Artin invariant
$\sigma_0$ determines $\NS(X)$ up to isometry, and we denote this lattice 
by $\Lambda_{p,\sigma_0}$.
Let $U'$ be the rank $2$ lattice with basis $\{Z,E\}$ and intersection matrix
$$
\left(
\begin{array}{cc}
  -2 & 1 \\
  1 & 0
\end{array}
\right) \,.
$$
To show the existence of a Jacobian elliptic fibration on $X$,
it suffices to find an isometric embedding of $U'$
into $\Lambda_{p,\sigma_0}$.
Since $U'$ is isometric to a hyperbolic plane $U$, and
since $\Lambda_{p,\sigma_0}$ is a sublattice of $\Lambda_{p,\sigma_0-1}$
for every $\sigma_0\geq2$, 
it suffices to show that $\Lambda_{p,9}$ contains $U$
in order to establish claim (1).
However, this follows from the explicit classification 
of the lattices $\Lambda_{p,\sigma_0}$ in \cite[Section 1]{rudakov shafarevich}:
namely, there exists an isometry
$$
   \Lambda_{p,9} \,\iso\, U\oplus H_p\oplus\left( I(-p)^{16} \right)_\ast,
$$
where the other lattices are defined and explained
in \cite[Section 1]{rudakov shafarevich}.
\qed\medskip

\begin{Remark}
  In characteristic $p\leq3$, we leave it to the reader to show the
  following if $X$ is a Shioda-supersingular K3 surface:
  \begin{enumerate}
     \item If $\sigma_0\leq9$, then $X$ admits a Jacobian genus $1$ fibration.
     \item If $\sigma_0=10$, then $X$ does not admit a Jacobian genus $1$ fibration.
       Moreover, if $p=3$ and $\sigma_0=6$, then $X$ does not admit a Jacobian
       quasi-elliptic fibration.
   \end{enumerate} 
 \end{Remark}

\subsection{Small Characteristics}
\label{subsec: small char}
Unfortunately, Proposition \ref{prop: spread out} rests on a theorem of Rudakov
and Shafarevich \cite{rudakov shafarevich degeneration} 
that supersingular K3 surfaces have potential good reduction,
which (currently) requires the assumption $p\geq5$.

\section{Moduli spaces}
\label{sec: ogus moduli}

In this section, we study the moving torsor families
from Proposition \ref{prop: spread out} using moduli spaces.
In order to avoid technical difficulties, we work with
moduli spaces of rigidified K3 crystals rather than 
moduli spaces of marked supersingular K3 surfaces.
As an interesting byproduct, we show that moduli spaces 
of rigidified K3 crystals are related
to each other by iterated $\PP^1$-bundles, together with a moduli 
interpretation.
In particular, this gives a new description of these moduli spaces, 
see Remark \ref{rem: crystal geometry}  and Remark  \ref{rem: twistor}.

\subsection{Recap of Ogus' period map}
In this subsection, we shortly review Ogus' articles
\cite{ogus} and \cite{ogus torelli}.
Let $N$ be a {\em supersingular K3 lattice}, that is,
the N\'eron--Severi lattice of a supersingular K3 surface
in characteristic $p$.
By \cite[Section 1]{rudakov shafarevich}, such a lattice
is determined up to isometry by $p$ and its Artin invariant $\sigma_0$.

\begin{Definition}
  Let $N$ be a supersingular K3 lattice.
  An {\em $N$-marked supersingular K3 surface} is
  a K3 surface $X$ together with an isometric embedding
  $N\to\NS(X)$.
\end{Definition}

We now assume $p\geq5$.
In \cite[Theorem (2.7)]{ogus torelli}, Ogus showed the existence of
a fine moduli scheme ${\cal S}_N$ for $N$-marked supersingular K3 surfaces,
and proved that it is locally of finite presentation, 
locally separated, and smooth of dimension
$\sigma_0(N)-1$ over $\FF_p$.
Moreover, ${\cal S}_N$ is almost proper, but neither of finite type nor
separated over $\FF_p$.
Here, we call a scheme {\em almost proper}, if it satisfies
the existence part of the valuative criterion for properness
with DVR's as test rings.

A {\em K3 crystal} of rank $22$ consists of
a triple $(H,\langle-,-\rangle,\Phi)$,
where $H$ is free $W$-module of rank $22$, $\langle-,-\rangle$ is
a symmetric bilinear form on $H$, and $\Phi$ is a Frobenius-linear 
endomorphism of $H$, that satisfies the conditions of
\cite[Definition 3.1]{ogus}.
For example, the $F$-crystal arising from $\Hcris{2}$ of a K3 surface, 
together with the symmetric bilinear
form coming from Poincar\'e duality, is a K3 crystal.
In case $H$ is of slope one, the K3 crystal is called
{\em supersingular}.
By the crystalline Torelli theorem \cite[Theorem I]{ogus torelli}, 
a supersingular K3 surface in characteristic $p\geq5$ 
is determined up to isomorphism 
by its supersingular K3 crystal.

In order to obtain Ogus' period map, we first have to rigidify
the K3 crystals:
by definition, the {\em Tate-module} of a K3 crystal $H$ is
defined to be $T_H:=\{x\in H:\Phi(x)=px\}$.
If $H$ is supersingular, then
$T_H$  is a free $\ZZ_p$-module of rank $22$, and the
bilinear form $\langle-,-\rangle$ on $H$ induces
a non-degenerate and non-perfect bilinear form on $T_H$.
Moreover, an $N$-marking of a supersingular K3 surface induces,
via the crystalline Chern map,
an isometric embedding of $N$ into the Tate-module of the 
associated K3 crystal, which motivates the following definition.

\begin{Definition}
  \label{def: rigidified crystal}
  Let $N$ be a supersingular K3 lattice.
  An {\em $N$-rigidified K3 crystal} is a pair $(\imath:N\to T_H, H)$, where
  $H$ is a K3 crystal, and $\imath$ is an isometric embedding.
\end{Definition}

By \cite[Proposition 4.6]{ogus},
there exists a moduli space ${\cal M}_N$
of $N$-rigidified K3 crystals, which is smooth and projective
of dimension $\sigma_0(N)-1$ over $\FF_p$.
We refer to Remark \ref{rem: crystal geometry} and the references given there
for details about its geometry.
Assigning to an $N$-marked supersingular K3 surface its $N$-rigidified
K3 crystal induces a morphism $\pi:{\cal S}_N\to{\cal M}_N$.

In order to get the period map, we have to equip $N$-rigidified
K3 crystals with ample cones, and refer to 
\cite[Definition 1.15]{ogus torelli} for definitions.
There exists a moduli scheme ${\cal P}_N$ of $N$-rigidified K3 crystals
with ample cones, which is
almost proper and locally of finite type over $\FF_p$.
Forgetting the ample cone induces an \'etale and surjective morphism
$f_N:{\cal P}_N\to{\cal M}_N$,
which is neither of finite type, nor separated.
Finally, assigning to an $N$-marked supersingular K3 surface
its $N$-rigidified supersingular K3 crystal together with the ample cone arising
from the ample cone of $X$ defines a lift of
$\pi$ to a morphism 
$$
   \widetilde{\pi}\,:\,{\cal S}_N\,\longrightarrow\,{\cal P}_N.
$$
This is Ogus' {\em period map}, and it is an isomorphism by \cite[Theorem III']{ogus torelli}.

\subsection{Moduli spaces of rigidified K3 crystals}
After these preparations, we now interpret Proposition \ref{prop: spread out}
in terms of rigidified K3 crystals:
if $X$ is a Jacobian elliptic fibration on a supersingular
K3 surface $X$, and $\oY\to C$ is as in Proposition \ref{prop: spread out}, 
then we obtain orthogonal decompositions 
$$ 
   \NS(X)\,\iso\, U\oplus \Lambda\mbox{ \qquad and \qquad }
   \NS(\oY_{\overline{\eta}})\,\iso\,U(p)\oplus\Lambda.
$$
More precisely, $U$ is the hyperbolic plane generated by the classes of a fiber $E$ 
and the zero-section $Z$ of the fibration, $\Lambda$ is defined to be $U^\perp$ 
inside $\NS(X)$, and  $U(p)$ is the lattice generated by $E$ and $pZ$. 
Then, we have the following theorem on moduli spaces of rigidified K3 crystals, 
which depends on these lattice decompositions only, and which is independent 
from Section \ref{sec:torsors}.
In Theorem \ref{thm: moduli interpretation} below, we will show that it is
indeed a manifestation of Proposition \ref{prop: spread out} 
on the level of K3 crystals.

\begin{Theorem}
   \label{thm: bundle structure}
   Let $N$ and $N_+$ be the supersingular 
   K3 lattices in odd characteristic $p$ 
   of Artin-invariants $\sigma_0$ and $\sigma_0+1$, respectively.
   Then, there exists a rank $20$ lattice $\Lambda$, and 
   orthogonal decompositions
   $$
       N\,\iso\, U\oplus \Lambda
       \mbox{ \qquad and \qquad }
       N_+\,\iso\, U(p)\oplus\Lambda\,,
   $$
   where $U$ denotes the hyperbolic plane.
   These decompositions gives rise to 
   a surjective morphism $\varpi_N$
   of moduli spaces of rigidified K3 crystals
   with a section $\sigma_N$
   $$
     \xymatrixcolsep{3pc}
     \xymatrix{
       {\cal M}_{N_+} \ar[d]^{\varpi_N} \\
       {\cal M}_N \ar@/^/[u]^{\sigma_N}  }
   $$
   which turns ${\cal M}_{N_+}$ into a $\PP^1$-bundle over ${\cal M}_N$.
 \end{Theorem}
 
\prf
We proceed in several steps.
\medskip

{\sc Step 1:} Setting up the lattices.

Since $\sigma_0\leq9$, there exists an isometry
$N\iso U\oplus\Lambda$ (see, for example, the proof of
Proposition \ref{prop:elliptic fibrations}).
Next, we choose a basis $\{E,Z\}$ of $U$ such that $E^2=0$, $Z^2=-2$, $E\cdot Z=1$.
Then, $E$ and $D:=pZ$ span a sublattice of $U$, which is isometric to $U(p)$.
Since $U(p)\oplus\Lambda$ is a supersingular K3 lattice of Artin invariant
$\sigma_0+1$, it is isometric to $N_+$ by the uniqueness result in
\cite[Section 1]{rudakov shafarevich}.
Thus, we obtain a commutative diagram of embeddings of lattices:
$$
\begin{array}{ccc}
    N_+ & \to & N \\
    \uparrow &  & \uparrow\\
    U(p) &\to & U
\end{array}
$$
\medskip

{\sc Step 2:} Translation from crystals into characteristic subspaces.

For the explicit computations, it is more convenient to work with characteristic subspaces 
rather than rigidified K3 crystals, and we
refer to \cite[Proposition 4.3]{ogus} for the translation between these two
points of view.
As in loc. cit., we define
$$
  N_0:=pN^\vee/pN
  \mbox{ \qquad and \qquad }
  (N_+)_0:=pN_+^\vee/pN_+,
$$
which are $\FF_p$-vector spaces of dimensions 
$2\sigma_0$ and $2\sigma_0+2$, respectively.
Also, the intersection forms turn $pN^\vee$ and $pN_+^\vee$ into 
sublattices of $N$ and $N_+$, respectively.
Moreover, by \cite[Proposition 3.13]{ogus},
the intersection forms on $N$ and $N_+$ are divisible by $p$
on $pN^\vee$ and $pN_+^\vee$, and induce perfect forms on $N_0$
and $(N_+)_0$.
A straight forward computation shows that  
the embedding $U(p)\subset N_+$ induces an isometry
$(N_+)_0\iso N_0\oplus (U\otimes\FF_p)$, where
$U\otimes\FF_p$ is generated by the classes of $D$ and $E$.
Tensoring the inclusion $N_+\subset N$ with $\FF_p$, 
we obtain a map $\gamma: N_+\otimes\FF_p\to N\otimes\FF_p$,
which has a one-dimensional kernel generated by $D$,
and whose cokernel is one-dimensional generated by $Z$.
Combining the remarks and computations of the previous
paragraph, we 
obtain a commutative diagram of $\FF_p$-vector spaces
$$
   \begin{array}{ccccc}
     (N_+)_0  & \iso & N_0 \oplus (U\otimes \FF_p) &\subset& N_+\otimes\FF_p \\
     & &  & & \downarrow \gamma \\
     & & N_0 &\subset& N\otimes\FF_p
   \end{array}
$$
For a field $k$ of characteristic $p$, we
set $\varphi := {\rm id}\otimes F_k^*$ on $N_0\otimes k$,
where $F_k$ denotes Frobenius.
By \cite[Definition 3.19]{ogus}, a 
{\em characteristic subspace} of $N_0\otimes k$ is a 
totally isotropic $k$-subvector space $K$, such that
$K+\varphi(K)$ is of dimension $\sigma_0+1$.
It is called {\em strictly characteristic} if it is characteristic and moreover 
$\sum_{i=0}^\infty \varphi^i(K)$ is equal to $N_0\otimes k$.
If $A$ is an $\FF_p$-algebra, then a {\em geneatrix} of $N_0\otimes A$ is
a direct summand of rank $\sigma_0$ such that the intersection form restricted
to it is identically zero, see \cite[p. 40]{ogus}.
Finally, a geneatrix is called {\em characteristic} if $K+F_A^*(K)$ is a direct summand
of rank $\sigma_0+1$.
Then, ${\cal M}_N$ represents the functor taking $A$ to the set of 
characteristic geneatrices of $N_0\otimes A$, see \cite[Proposition 4.6]{ogus}.
\medskip

{\sc Step 3:} Definition of $\sigma_N$.

By \cite[Definition 4.1]{ogus}, 
${\cal M}_N$ parametrizes $N$-rigidified K3 crystals, that is, 
pairs $(\imath: N\to T_H, H)$ as in Definition \ref{def: rigidified crystal}.
Composing $\imath$ with $N_+\to N$ turns an $N$-rigidified
K3 crystal into an $N_+$-rigidified K3 crystal, which defines a 
morphism $\sigma_N:{\cal M}_N\to{\cal M}_{N_+}$.
Translated into geneatrices, this becomes the following:
if $A$ is an $\FF_p$-algebra, and if $K\subset N_0\otimes A$ 
is a characteristic geneatrix, 
then $\gamma^{-1}(K)\iso K\oplus (D\cdot A)$,
and easily seen to be a characteristic geneatrix of $(N_+)_0\otimes A$.
Using \cite[Proposition 4.3]{ogus}, it is not difficult to see
that the assignment
$$
    K\,\mapsto\,\gamma^{-1}(K)
$$
describes $\sigma_N$ in terms of characteristic geneatrices.
\medskip

{\sc Step 4:} Definition of $\varpi_N$.

For an $\FF_p$-algebra $A$, and a direct summand 
$K\subseteq(N_+)_0\otimes A$, we set
$$
  \Gamma_+(K) \,:=\, {\rm pr}_{N_0}\left(K\cap (E^\perp\otimes A)\right),
$$
where ${\rm pr}_{N_0}$ denotes the projection $(N_+)_0\otimes A\to N_0\otimes A$.
A straight forward calculation shows that if $K$ is a characteristic geneatrix
of $(N_+)_0\otimes A$, then $\Gamma_+(K)$ is a characteristic geneatrix 
of $N_0\otimes A$.
Thus, the assignment
$$
  K\,\mapsto\,\Gamma_+(K)
$$
defines a morphism ${\cal M}_{N_+}\to{\cal M}_N$ that we denote by
$\varpi_N$.
\medskip

{\sc Step 5:} $\sigma_N$ is a section of $\varpi_N$.

If $A$ is an  $\FF_p$-algebra and $K$ is a direct summand 
of  $N_0\otimes A$, then it follows from the
definitions that  $\Gamma_+(\gamma^{-1}(K))=K$,
which shows that $\varpi_N\circ\sigma_N={\rm id}$.
\medskip

{\sc Step 6:} $\varpi_N$ defines a $\PP^1$-bundle structure.

Using the isomorphism
$(N_+)_0\iso N_0\oplus (U\otimes\FF_p)$,
we have a projection
$$
  {\rm pr}_U\,:\,
  (N_+)_0\,\to\, (U\otimes \FF_p)\,.
$$
Now, let $k$ be an algebraically closed field of characteristic $p$,
and let  $K_0\subset N_0\otimes k$ be a characteristic subspace, that is, 
a $k$-rational point of ${\cal M}_N$.
A straight forward computation shows that if
$K_+\subset (N_+)_0\otimes k$ is characteristic,
then ${\rm pr}_U(K_+\cap\varphi(K_+))$ is one-dimensional.
This shows that
$K_+\cap\varphi(K_+)\cap (N_0\otimes k)$
is $(\sigma_0-1)$-dimensional, where we view
$N_0$ again as a subspace of $(N_+)_0$.
In particular, if $\Gamma_+(K_+)=K_0$, 
then $K_+\cap\varphi(K_+)\cap (N_0\otimes k)=K_0\cap\varphi(K_0)$.
Thus, every characteristic subspace $K_+\subset(N_+)_0\otimes k$
with $\Gamma_+(K_+)=K_0$ contains the $(\sigma_0-1)$-dimensional 
and totally isotropic subspace $K_0\cap\varphi(K_0)$.
Let $k_1,...k_{\sigma_0-1}$ be a basis of $K_0\cap\varphi(K_0)$,
and choose $v\in K_0$ such that $K_0=\langle v, K_0\cap\varphi(K_0)\rangle$
and $\varphi(K_0)=\langle \varphi(v), K_0\cap\varphi(K_0)\rangle$.
We normalize $v$ such that $\langle v, \varphi(v)\rangle=1$.
Then, another straight forward calculation shows that 
$K_+\subset(N_+)_0\otimes k$ is characteristic with $\Gamma_+(K_+)=K_0$ 
if and only if either $K_+=\langle K_0, E\rangle$ or if
there exists a unique $\lambda\in k$ such that
$$
 K_+ \,=\, \left\langle\, 
 k_1,\,...,\,k_{\sigma_0-1},\, v+\lambda E,\, v-\lambda \varphi(v)+D+\lambda E\, 
 \right \rangle \,.
$$
Thus, the fiber of $\varpi_N$ over $K_0$ is isomorphic to $\PP^1$,
and since $K_0$ was chosen arbitrarily, this shows that all 
fibers of $\varpi_N$ over geometric points of ${\cal M}_N$
are isomorphic to $\PP^1$.
In particular, $\varpi_N$ is a conic bundle.
Since $\sigma_N$ is a section of $\varpi_N$, this conic bundle is 
a $\PP^1$-bundle.
\qed\medskip
 
\begin{Remark}
  \label{rem: crystal geometry}
  In \cite[Examples 4.7]{ogus}, Ogus explicitly described
  ${\cal M}_N$ in the following  cases
  $$
     \begin{array}{cc}
        \sigma_0(N) & {\cal M}_N\\
        \hline
      1 & \Spec \FF_{p^2}\\
      2 & \PP^1 \,\times\, \Spec\FF_{p^2}\\
      3 & (\PP^1\times\PP^1)\,\times\, \Spec\FF_{p^2}  
     \end{array}
  $$
  By our previous theorem, ${\cal M}_N$ is an iterated $\PP^1$-bundle
  over $\Spec\FF_{p^2}$, and we refer to
  \cite[Remark 4.8]{ogus} and \cite[Theorem 3.21]{ogus} 
  for further descriptions.
\end{Remark}

\subsection{The moduli interpretation}
The previous theorem is about moduli spaces of
rigidified K3 crystals.
The following theorem links it to the moving torsor families of 
supersingular K3 surfaces from Proposition \ref{prop: spread out},
and gives a moduli interpretation of $\varpi_N$ and $\sigma_N$.

\begin{Theorem}
   \label{thm: moduli interpretation}
   We keep the notations and assumptions of Theorem \ref{thm: bundle structure}.
   Moreover, we assume $p\geq5$.
   \begin{enumerate}
    \item Let $X$ be a supersingular K3 surface with $\NS(X)\iso N$, and let
      $[X]\in{\cal M}_N$ be the associated K3 crystal.
      Then, the family 
      $$
          \oY\,\to\, C,
      $$
      from Proposition \ref{prop: spread out} comes with an $N_+$-marking,
      such that the associated family of $N_+$-rigidified K3 crystals maps 
      onto $\varpi_N^{-1}([X])$.
      \item Being the fiber over $0\in C$, the surface $X$ inherits an $N_+$-marking,
        and the corresponding K3 crystal is $\sigma_N([X])$.
   \end{enumerate}
\end{Theorem}

\prf
We keep the notations from the proof of Theorem \ref{thm: bundle structure}.
Given $X$ as in claim (1), we choose the isomorphism $N\iso\NS(X)$ such that
$U\subset N$ corresponds to a Jacobian elliptic fibration on $X$,
see also the proof of Proposition \ref{prop:elliptic fibrations}.
Next, let $\oY\to C$ be the associated family 
from Proposition \ref{prop: spread out}.
Let $\eta\in C$ be the generic point, set $R:=\OO_{C,0}$, 
choose a uniformizer $t\in R$, and note that $k(C)$
is the field of fractions of $R$.
By Proposition \ref{prop:picard index p},
the isomorphism $N\iso\NS(X)$ induces an isomorphism
$N_+\iso\NS(\oY_{\eta})$, and, via restriction, 
the whole family $\oY\to C$ becomes $N_+$-marked.
More precisely, we obtain orthogonal decompositions
of $N$ and $N_+$  and an embedding $N_+\subset N$ 
as in step 1 of the proof of Theorem \ref{thm: bundle structure}.

Let us now determine the characteristic subspaces associated
to $X$ and $\oY_\eta$.
As explained on \cite[p. 365]{ogus torelli}, these 
arise as kernels of the Chern class $c_{\rm dR}$.
We have a commutative diagram
$$
\begin{array}{ccccc}
  \NS(\oY_{\eta}) &\to&  \NS(\oY_{\eta})\otimes_\ZZ k(C) &\stackrel{c_{\rm dR}}{\longrightarrow}& \HdR{2}( \oY_{\eta}/\eta) \\
  \uparrow && \uparrow && \uparrow \\
  \NS(\oY_{R}) &\to&  \NS(\oY_{R})\otimes_\ZZ R &\stackrel{c_{\rm dR}}{\longrightarrow}& \HdR{2}( \oY_{R}/R) \\
  \downarrow && \downarrow{\scriptstyle\gamma'} && \downarrow \\
  \NS(X) &\to &\NS(X)\otimes_\ZZ k & \stackrel{c_{\rm dR}}{\longrightarrow}& \HdR{2}(X/k)
\end{array}
$$
whose vertical arrows are restriction maps.
Let $K_0':=\varphi^{-1}(K_0)\subset N_0\otimes k$ be the characteristic subspace 
associated to $X$.
It is not difficult to see that there exists a lift of 
$K_0'$ to an $R$-submodule $\widetilde{K}_0'\subset N_+\otimes R$ of rank $\sigma_0$
that is contained in $\ker(c_{{\rm dR}})$.
More precisely, if $k_1,...,k_{\sigma_0}$ is a basis of $K_0'$, and
$\overline{k}_i:=k_i\otimes 1\in N_0\otimes R$, there exist lifts of the $k_i$ to
$\ker(c_{\rm dR})$ of the form
$$
 \overline{k}_i + t \overline{n}_i + \alpha_i D + t\beta_i E,\mbox{ \qquad } i=1,...,\sigma_0,
$$
where $\overline{n}_i\in N_0\otimes R$, and $\alpha_i, \beta_i\in R$.
There is one more element in $\ker(c_{\rm dR})$, linearly independent from these, and
without loss of generality, it is not divisible by $t$ and lies in the kernel of $\gamma'$.
Thus, we may choose it to be of the form
$$
   t\overline{n}_0+D+t\beta E,
$$
where $\overline{n}_0\in N_0\otimes R$ and $\beta\in R$.
Since these $\sigma_0+1$ elements lie inside $\ker(c_{dR})$, they
form a totally isotropic subspace.
After some tedious computations exploiting this isotropy, we find
that $\ker(c_{dR})$ contains a free $R$-submodule 
$\widetilde{K}_+'$ of rank $\sigma_0+1$ generated by elements of the form
$$
\begin{array}{ccccc}
  \overline{k}_i    &  &     &+&  t\mu_i E \\
  t\overline{n}_0 &+& D &+& t\beta E
\end{array}
$$
Simply for dimensional reasons,
$\widetilde{K}_+'\otimes\overline{k(C)} \subset (N_+)_0\otimes \overline{k(C)}$ is
the characteristic subspace associated to $\oY_{\overline{\eta}}$.

Using this explicit description, we compute $\gamma'(\widetilde{K}_+')=K_0'$
and $\Gamma_+(\widetilde{K}_+')=K_0'\otimes_k R$,
where $\Gamma_+$ is defined as in the proof of step 4 of 
Theorem \ref{thm: bundle structure}.
In particular, the classifying map $f_C:C\to{\cal M}_{N_+}$ 
maps to the fiber $\varpi_N^{-1}([X])$.
Since the fibers of $\varpi$ and $C$ are proper irreducible curves
and $f_C$ is not constant, $f_C$ maps
surjectively onto $\varpi_N^{-1}([X])$, which establishes claim (1).
The fiber over $0\in C$ is isomorphic to $X$,
and the $N_+$-marking of $\NS(X)$ induced from the $N_+$-marking of the 
family $\oY\to C$ arises via $N_+\to N\iso\NS(X)$.
Thus, by the definition of $\sigma_N$ in step 3 of the proof of
Theorem \ref{thm: bundle structure}, the associated 
$N_+$-rigidified K3 crystal is $\sigma_N([X])$, which establishes claim (2).
\qed\medskip

Let us give an immediate corollary of
Theorem \ref{thm: bundle structure} and Theorem \ref{thm: moduli interpretation},
which is actually everything we will need to prove the results of the next section.

\begin{Corollary} 
     \label{cor: ogus}
    Let $Y$ be a supersingular K3 surface in characteristic $p\geq5$
    with $\sigma_0(Y)\geq2$.
    Then, there exists a supersingular K3 surface $X$ with 
    $\sigma_0(X)=\sigma_0(Y)-1$ that is purely inseparably isogenous
    of height $2$ to $Y$.
\end{Corollary}

\prf
By Theorem \ref{thm: bundle structure} and Theorem \ref{thm: moduli interpretation},
there exists a family $\oY\to C$ of $N_+$-marked supersingular K3 surfaces,
where $\sigma_0(N_+)=\sigma_0(Y)$, such that
$Y$ is a member of this family.
Let $X$ be the special fiber over $0\in C$, and then, the assertion follows
from Proposition \ref{prop: spread out}.
\qed\medskip

\begin{Remark}
  \label{rem: twistor}
  The unirationality of ${\cal M}_N$ is clear from 
  Ogus' description \cite[Theorem 3.21]{ogus}, whereas our description as iterated
  $\PP^1$-bundle is new.
  It is likely that Theorem \ref{thm: bundle structure} and Theorem \ref{thm: moduli interpretation}
  extend in some form to the moduli spaces ${\cal S}_N$ of $N$-marked supersingular K3 surfaces.
  However, since these latter spaces are neither of finite type nor separated, the proofs and maybe even
  the statements would probably be rather technical and involved.
  Much better behaved are moduli spaces of polarized K3 surfaces.
  In \cite[Section 9]{lieblich ubiquity}, Lieblich announced the existence of families
  of supersingular K3 surfaces over $\Aff^1$ using moduli spaces of twisted sheaves.
  As an application, he announces the uniruledness of the supersingular loci
  of moduli spaces of polarized K3 surfaces.
  For complex K3 surfaces, non-trivial families over $\PP^1$, whose general member
  is not algebraic, arise from twistor spaces, see \cite[Section 25]{huybrechts}.
  For example, Markman \cite[Section 7]{markman} and 
  Verbitsky \cite[Section 1.4]{verbitsky} studied
  twistor spaces together with Lagrangian fibrations, which  
  is similar to our moving torsor families.
\end{Remark}

\subsection{Small Characteristics}
\label{subsec: ogus small char}
The results of this section build on Ogus' articles
\cite{ogus} and \cite{ogus torelli}.
In \cite{ogus}, he develops the theory of supersingular K3 crystals,
and the assumption $p\geq3$ is built in from
the very beginning: quadratic and symplectic forms 
play an important role, which is why characteristic $2$ is excluded.
In \cite{ogus torelli}, $p\geq5$ had to be assumed, not only
because it rests on \cite{ogus}, but also since it needs the theorem of
Rudakov--Shafarevich \cite{rudakov shafarevich degeneration} 
on potential good reduction of supersingular K3 surfaces,
see \cite[p. 364]{ogus torelli}.

\section{Supersingular K3 surfaces are unirational}
\label{sec:shioda}

In this section, we prove that supersingular K3 surfaces 
in characteristic $p\geq5$ are related by
purely inseparable isogenies, which is
an analog of the Shioda--Inose structure theorem for singular K3 surfaces,
see Theorem \ref{shioda inose in p}.
Since Shioda \cite{shioda some results}
showed that supersingular Kummer surfaces are unirational, we 
deduce the Artin--Rudakov--Shafarevich--Shioda
conjecture on unirationality of all supersingular
K3 surfaces. 
Finally, we treat unirationality of Enriques surfaces.

\subsection{Isogenies between supersingular K3 surfaces}
We now come to the main theorem of this article,
which is a structure result for supersingular K3 surfaces.
We note that Rudakov and Shafarevich conjectured this
already in Question 8 at the end of \cite{rudakov shafarevich}.
We refer to Section \ref{subsec: isogeny} for the connection with
a conjecture of Shafarevich about isogenies between complex K3 surfaces.

\begin{Theorem}
 \label{thm:correspondences}
  Let $X$ and $X'$ be supersingular K3 surfaces with Artin invariants $\sigma_0$ and
  $\sigma_0'$ in characteristic $p\geq5$.
  \begin{enumerate}
    \item There exist purely inseparable isogenies
      $$
        X \,\dashrightarrow\,X'\,\dashrightarrow\,X\,,
      $$
      both of which are of height $2\sigma_0+2\sigma_0'-4$.
   \item Let $E$ be a supersingular elliptic curve. Then, there exist isogenies
     $$
     {\rm Km}(E\times E)\,\dashrightarrow\,X\,\dashrightarrow\,{\rm Km}(E\times E)\,,
     $$ 
     both of which are purely inseparable of height $2\sigma_0-2$.
  \end{enumerate}
\end{Theorem}

\prf
If $\sigma_0\geq2$, then
there exists a supersingular K3 surface with Artin invariant $\sigma_0-1$
that is purely inseparable isogenous of height $2$ to $X$
by Corollary \ref{cor: ogus}.
By induction, we obtain a purely inseparable isogeny $\varphi$
of height $2\sigma_0-2$ from $X$ to a supersingular
K3 surface with Artin invariant $\sigma_0=1$.
However, there exists only one such surface, namely the Kummer surface
${\rm Km}(E\times E)$, where 
$E$ is a supersingular elliptic curve \cite[Corollary (7.14)]{ogus}.
Since the $(2\sigma_0-2)$-fold Frobenius of $X$ factors through $\varphi$,
we obtain claim (2).

By the established claim (2), there exists a purely inseparable isogeny 
$\varphi':{\rm Km}(E\times E)\dashrightarrow X'$ of height $2\sigma_0'-2$.
Then, $\varphi'\circ\varphi$ is a purely inseparable isogeny 
$X\dashrightarrow X'$ of height $2\sigma_0+2\sigma_0'-4$.
As before, the $(2\sigma_0+2\sigma_0'-4)$-fold Frobenius of $X$ 
factors through $\varphi'\circ\varphi$ and we obtain claim (1).
\qed\medskip

\begin{Remark}
  \label{rem:correspondences}
  Naively, one might expect that K3 surfaces of Picard rank $\geq\rho$ form 
  a codimension $\rho$ subset inside the 
  moduli space.
  This expectation is fulfilled for 
  singular K3 surfaces ($\rho=20$), since they form a countable set.
  But then, one would expect that K3 surfaces with $\rho=22$ should
  not exist at all, and the fact that they come in $9$-dimensional families
  is even more puzzling.
  However, by Theorem \ref{thm:correspondences}, there exists only one
  supersingular K3 surface in every positive characteristic up to purely inseparable
  isogeny.
  By Proposition \ref{prop: spread out}, these isogenies come in families,
  which gives an explanation why supersingular K3 surfaces form 
  $9$-dimensional moduli spaces.
 \end{Remark}

\subsection{Supersingular K3 surfaces are unirational}
Since Shioda \cite{shioda some results} 
showed that supersingular Kummer
surfaces are unirational, the previous theorem  implies 
the conjecture of Artin, Rudakov, Shafarevich, and Shioda.

\begin{Theorem}
 \label{thm:main}
  Supersingular K3 surfaces in characteristic $p\geq5$ are unirational.
\end{Theorem}

\prf
In odd characteristic, supersingular 
Kummer surfaces are unirational by \cite[Theorem 1.1]{shioda some results}.
The assertion then follows from Theorem \ref{thm:correspondences}.
\qed\medskip

We recall that a surface is called a {\em Zariski surface} if there exists a
dominant, rational, and purely inseparable map of degree $p$ from $\PP^2$ onto it.
Although the map from $\PP^2$ onto a supersingular Kummer surface
constructed by Shioda in \cite{shioda some results} is inseparable, 
it is not purely inseparable.
Using a different construction, Katsura \cite[Theorem 5.10]{katsura} showed that 
supersingular Kummer surfaces with $\sigma_0=1$ in characteristic 
$p\not\equiv1\mod 12$ are Zariski surfaces.
This strengthens Theorem \ref{thm:main}, and gives
a partial answer to a question of Rudakov and Shafarevich, who asked 
and actually doubted whether supersingular K3 surfaces are 
purely inseparably unirational, see
Question 6 at the end of \cite{rudakov shafarevich}.

\begin{Corollary}
  A supersingular K3 surface in characteristic $p\geq5$ with $p\not\equiv1\mod12$
  is purely inseparably unirational.\qed
\end{Corollary}

In Section \ref{subsec: brauer}, we discussed different notions of supersingularity
for K3 surfaces and the relation to the Tate-conjecture.
Combining Theorem \ref{thm: equivalence} and 
Theorem \ref{thm:main}, we obtain the following equivalence.

\begin{Theorem}
 \label{thm: conclusion}
 For a K3 surface $X$ in characteristic $p\geq5$, the following
 conditions are equivalent:
 \begin{enumerate}
  \item $X$ is unirational.
  \item The Picard rank of $X$ is $22$.
  \item The formal Brauer group of $X$ is of infinite height.
  \item For all $i$, the $F$-crystal $\Hcris{i}(X/W)$ is of slope $i/2$.
 \end{enumerate}
\end{Theorem}

\prf
If $X$ is unirational, then its Picard rank is $22$ by \cite[Corollary 2]{shioda an example},
which establishes (1)$\Rightarrow$(2).
The converse direction (2)$\Rightarrow$(1) is Theorem \ref{thm:main}.
The equivalences (2)$\Leftrightarrow$(3)$\Leftrightarrow$(4)
are Theorem \ref{thm: equivalence}.
\qed\medskip

\subsection{Enriques surfaces}
As a consequence of Theorem \ref{thm:main}, we now characterize
the unirational ones among Enriques surfaces, which generalizes a result of 
Shioda \cite[Theorem 3.3]{shioda some results}.

\begin{Theorem}
   An Enriques surface $X$ in characteristic $p\geq2$ is unirational if and only if
   \begin{enumerate}
     \item $p=2$ and $X$ is not singular (that is, $\Pic^\tau_{X/k}\neq\mu_2$), or
     \item $p\neq2$ and the covering K3 surface is supersingular.
   \end{enumerate}
\end{Theorem}

\prf
Assertion (1) is shown in \cite[Corollary I.1.3.1]{cd}.

By \cite[Lemma 3.1]{shioda some results}, an Enriques surface $X$
in characteristic $p\geq3$ is unirational if and only if its covering 
K3 surface $\widetilde{X}$ is unirational.
Thus, if $p\geq5$, then assertion (2) follows from Theorem \ref{thm: conclusion}.
If $p=3$ and $X$ is unirational, then $\widetilde{X}$ is unirational,
and thus, supersingular.
Conversely, if $p=3$ and $\widetilde{X}$ is supersingular, then
$\sigma_0(\widetilde{X})\leq5$ by \cite[Corollary 3.4]{jang}
and thus, $\widetilde{X}$ is unirational by \cite{rudakov shafarevich},
which implies the unirationality of $X$.
\qed\medskip

\subsection{Small Characteristics}
\label{subsec: small char 2}
As in Section \ref{subsec: small char} and Section \ref{subsec: ogus small char},
let us discuss what we know and do not know in characteristic $p\leq3$.
  \begin{enumerate}
  \item Using quasi-elliptic fibrations, Rudakov and Shafarevich \cite{rudakov shafarevich}
      showed that Shioda-supersingular K3 surfaces in characteristic $2$
     and supersingular K3 surfaces with $\sigma_0\leq6$ in characteristic $3$
     are Zariski surfaces, and thus, unirational.
     Therefore,  the question remains
     whether supersingular K3 surfaces with $\sigma_0\geq7$ in characteristic $3$ 
     are unirational.
     By Proposition \ref{prop: deformation} together with the comments made in
     Section \ref{subsec: small char}, there exists at least
     a $6$-dimensional family
     of unirational K3 surfaces with $\sigma_0=7$ in characteristic $3$.
   \item Theorem \ref{thm:correspondences} rests on Corollary \ref{cor: ogus},
     and we refer to Section \ref{subsec: ogus small char} for details.
     On the other hand, quasi-elliptic K3 surfaces 
     are Zariski surfaces, and thus, related by purely inseparable isogenies.
   \item The implication (1)$\Rightarrow$(2) of Theorem \ref{thm: conclusion}
     holds in any characteristic and we  discussed it converse above.
     The implication (2)$\Rightarrow$(3) holds in any characteristic, and its
     converse would follow from the Tate-conjecture for K3 surfaces, which
     is true in characteristic $3$ by \cite{madapusi}. 
     The equivalence (3)$\Leftrightarrow$(4) holds in every characteristic.
\end{enumerate} 
In particular (see also Section \ref{subsec: ogus small char}), 
once supersingular K3 surfaces in characteristic $3$ are shown to have potential good reduction,
the results of this section will also hold in characteristic $3$.

\appendix
\section{Erratum}

As was pointed out by Bragg and Lieblich in \cite{BL}, the proof of
Proposition 3.5 contains a mistake, which has the effect that the main result,
the unirationality of supersingular K3 surfaces Theorem 5.3, 
remains a conjecture. 

\subsection*{Degree p multisections}
More precisely, the proof only shows the following:

\begin{Proposition}
  We keep the notations and assumptions of Proposition 3.2
  and assume that $\widehat{{\rm Br}}(X)$ is not $p$-divisible.
  Set $R:=k[[t]]$ and $S:=\Spec R$, and let
  $$
     \overline{A} \,\to\, Y\times_{k} S \,\to\, S
  $$ 
  be a family of $A$-torsors
  associated to a non-trivial $p$-torsion element of $\widehat{{\rm Br}}(X)(R)$
  as in  Proposition 3.2.
  Then, after possibly replacing $S$ by some finite flat cover
  \begin{enumerate}     
    \item There exists a degree $p$ multisection 
      $D_\xi\subset\overline{A}_\xi$
      such that the induced morphism
 $$
         D_\xi\,\to\, \xi\times_{k} S
$$
      is finite and flat of degree $p$.
    \item More precisely, $D_\xi\to\xi\times_k S$
      arises from a family of $A_\xi[p]$-torsors.
  \end{enumerate}
\end{Proposition}

In particular, contrary to what was claimed in  Proposition 3.5, one 
may not be able to choose the degree $p$ multisection $D_\xi$ to be radicial,
that is, purely inseparable, over $\xi\times_{k}S$.

The strategy of establishing  Theorem 5.3 consisted in showing
that any two supersingular K3 surfaces in characteristic $p\geq5$ can 
be connected by one-dimensional families of the type as constructed in Section 3
and then, that unirationality is preserved in these families.
This reduced the unirationality conjecture to the showing that there exists at least
one unirational K3 surface in each characteristic, where one can use Kummer surfaces,
whose unirationality was established by Shioda.

Now, in order to show that unirationality is preserved in these
one-dimensional families, one needs Assertion (4) of Theorem 3.6,
whose proof requires that one can find degree $p$ multisections 
in Proposition 3.5  that are radicial over the base.
Without this, one can only conclude that any two supersingular K3 surfaces
in characteristic $p\geq5$ are connected by a sequence of very special
degree $p$ correspondences (as in Assertion (3) of  Theorem 3.6),
which is not sufficient for unirationality.

\subsection*{Fibration structure}
Moreover, the proof of  Theorem 4.3  contains a mistake. 
In fact, the proof only shows the following:

\begin{Theorem}
   Let $N$ and $N_+$ be the supersingular 
   K3 lattices in odd characteristic $p$ 
   of Artin-invariants $\sigma_0$ and $\sigma_0+1$, respectively.
   Then, there exists a rank $20$ lattice $\Lambda$, and isometries
   $$
       N\,\iso\, U\oplus \Lambda
       \mbox{ \qquad and \qquad }
       N_+\,\iso\, U(p)\oplus\Lambda\,,
   $$
   where $U$ denotes the hyperbolic plane.
   This decomposition gives rise to 
   a dominant rational map $\varpi_N$
   of moduli spaces of rigidified K3 crystals
   with a section $\sigma_N$
   $$
     \xymatrixcolsep{3pc}
     \xymatrix{
       {\cal M}_{N_+} \ar[d]^{\varpi_N} \\
       {\cal M}_N \ar@/^/[u]^{\sigma_N}  }.
   $$
   More precisely,
   \begin{enumerate}
   \item this map is defined on the open set of ${\cal M}_{N+1}$ of rigidified K3 crystals
   with Artin invariant $\sigma_0+1$ and maps surjectively onto the open set
   of rigidified K3 crystals of Artin invariant $\sigma_0$.
   \item The normalization of the geometric generic fibre is $\PP^1$.
    (In general, $\varpi_N$ may not be generically smooth.)
   \end{enumerate}
 \end{Theorem}

In the proof, the map is only defined on geometric points and extension
to a morphism would require more arguments.
In particular, the computations in Step 6 of the proof only show that the normalisation
of the geometric generic fibre is isomorphic to $\PP^1$.
 
In particular, the table Remark 4.4 has to be corrected as follows:
for $\sigma_0(N)=3$, the moduli space ${\cal M}_N$ consists 
geometrically of two copies of the Fermat hypersurface of degree $(p+1)$ in $\PP^3$.
In this case, $\varpi_N$ arises geometrically by choosing a pencil of planes through a line of this
surface and then, one obtains a fibration over $\PP^1$ , whose geometric generic
fibre is a singular rational curve that is not smooth.
This gives an example, where the geometric generic fibre of $\varpi_N$ is not 
smooth, and in particular, not isomorphic to $\PP^1$.

\end{document}